\message{by Udo Hertrich-Jeromin, July 1997}\immediate\write1{}

\magnification = \magstep0
\overfullrule = 2pt
\vsize = 525dd
\hsize = 27.0cc
\topskip = 13dd
\hoffset = 1.8cm
\voffset = 1.8cm
\parindent = 0pt
\parskip = 1ex plus 3pt

\message{more fonts: petit and fraktur;}

\font\bfbig = cmbx10 scaled \magstep2   

\font\eightrm = cmr10 scaled 800        
\font\sixrm = cmr7 scaled 850
\font\fiverm = cmr5
\font\eighti = cmmi10 scaled 800
\font\sixi = cmmi7 scaled 850
\font\fivei = cmmi5
\font\eightit = cmti10 scaled 800
\font\eightsy = cmsy10 scaled 800
\font\sixsy = cmsy7 scaled 850
\font\fivesy = cmsy5
\font\eightsl = cmsl10 scaled 800
\font\eighttt = cmtt10 scaled 800
\font\eightbf = cmbx10 scaled 800
\font\sixbf = cmbx7 scaled 850
\font\fivebf = cmbx5

\font\fivefk = eufm5                    
\font\sixfk = eufm7 scaled 850
\font\sevenfk = eufm7
\font\eightfk = eufm10 scaled 800
\font\tenfk = eufm10

\newfam\fkfam
        \textfont\fkfam=\tenfk \scriptfont\fkfam=\sevenfk
                \scriptscriptfont\fkfam=\fivefk

\def\eightpoint{%
        \textfont0=\eightrm \scriptfont0=\sixrm \scriptscriptfont0=\fiverm
                \def\rm{\fam0\eightrm}%
        \textfont1=\eighti  \scriptfont1=\sixi  \scriptscriptfont1=\fivei
                \def\oldstyle{\fam1\eighti}%
        \textfont2=\eightsy \scriptfont2=\sixsy \scriptscriptfont2=\fivesy
        \textfont\itfam=\eightit \def\it{\fam\itfam\eightit}%
        \textfont\slfam=\eightsl \def\sl{\fam\slfam\eightsl}%
        \textfont\ttfam=\eighttt \def\tt{\fam\ttfam\eighttt}%
        \textfont\bffam=\eightbf \scriptfont\bffam=\sixbf
                \scriptscriptfont\bffam=\fivebf \def\bf{\fam\bffam\eightbf}%
        \textfont\fkfam=\eightfk \scriptfont\fkfam=\sixfk
                \scriptscriptfont\fkfam=\fivefk
        \rm}
\skewchar\eighti='177\skewchar\sixi='177\skewchar\eightsy='60\skewchar\sixsy='60

\def\fk{\fam\fkfam}
\def\petit{\eightpoint}%


\def\today{\ifcase\month\or
        January\or February\or March\or April\or May\or June\or
        July\or August\or September\or October\or November\or December\fi
        \space\number\day, \number\year}
\def\newline{\hfil\break}

\def\framebox#1{\vbox{\hrule\hbox{\vrule\hskip1pt
        \vbox{\vskip1pt\relax#1\vskip1pt}\hskip1pt\vrule}\hrule}}

\message{ lay out: title page,}

\def\Date#1{}

\newdimen\fullhsize \fullhsize=\hsize

\newbox\TITLEBOX \setbox\TITLEBOX=\vbox{}
\newdimen\TITLEHEIGHT
\long\def\Title#1#2#3#4{\setbox\TITLEBOX=\vbox{%
        \topglue3truecm%
        \noindent{\bfbig#1}\vskip12pt%
        \noindent{\bf#2}\vskip6pt%
        \noindent{\petit#3}\vskip10pt%
        \noindent{\petit\today}\vskip32pt%
        \noindent{\bf Summary.}\enspace#4\vskip32pt
        \relax}%
        \TITLEHEIGHT=\ht\TITLEBOX%
        \if\IFTHANKS Y\advance\TITLEHEIGHT by\ht\THANKBOX\fi%
        \advance\vsize by-\TITLEHEIGHT%
        \relax}

\newtoks\RunAuthor\RunAuthor={} \newtoks\RunTitle\RunTitle={}
\def\ShortTitle#1#2{\RunTitle={#1}\RunAuthor={#2}}
\headline={\ifnum\pageno=1 {\hfil} \else
        \ifodd\pageno {\petit{\the\RunTitle}\hfil\folio}
        \else {\petit\folio\hfil{\the\RunAuthor}} \fi \fi}
\def\makeheadline{\vbox{%
        \hbox to\fullhsize{\the\headline}%
        \vss\nointerlineskip\kern2pt%
        \hbox to\fullhsize{\hrulefill}\kern7pt}}

\newbox\THANKBOX \let\IFTHANKS N
\setbox\THANKBOX=\vbox{\kern12pt\hbox to\fullhsize{\hrulefill}\kern4pt}
\def\Thanks#1#2{\nobreak${}^{#1}$\global\let\IFTHANKS Y%
        \global\setbox\THANKBOX=\vbox{\parindent20pt\baselineskip9pt%
        \unvbox\THANKBOX{\petit\item{${}^{#1}$}{#2}} }}
\footline={\ifnum\pageno=1 \if\IFTHANKS Y\box\THANKBOX\fi
        \else\hfill\fi}
\def\makefootline{\hbox to\fullhsize{\the\footline}}


\def\OneColumn{\output={%
        \shipout\vbox{\makeheadline%
                \ifnum\pageno=1\hbox to\fullhsize{\box\TITLEBOX\hfil}\fi%
        \hbox to\fullhsize{\pagebody}%
        \makefootline}%
                \ifnum\pageno=1\global\advance\vsize by\TITLEHEIGHT\fi%
        \advancepageno%
        \ifnum\outputpenalty>-20000%
        \else\dosupereject%
        \fi}%
        }%
\OneColumn

\newbox\LEFTCOLUMN
\def\TwoColumn{\let\lr=L \hsize=.485\fullhsize \output={%
        \if L\lr\global\setbox\LEFTCOLUMN=\leftline{\pagebody}\global\let\lr=R%
        \else \shipout\vbox{\makeheadline%
                \ifnum\pageno=1\hbox to\fullhsize{\box\TITLEBOX\hfil}\fi%
        \hbox to\fullhsize{\box\LEFTCOLUMN\hfil\leftline{\pagebody}}%
        \makefootline}%
                \ifnum\pageno=1\global\advance\vsize by\TITLEHEIGHT\fi%
        \advancepageno\global\let\lr=L%
        \fi%
        \ifnum\outputpenalty>-20000%
        \else\dosupereject%
        \fi}%
        }%


\def\PLabel#1{\xdef#1{\nobreak(p.\the\pageno)}}

\message{headings and}

\newcount\SECNO \SECNO=0
\newcount\SUBSECNO \SUBSECNO=0
\newcount\SUBSUBSECNO \SUBSUBSECNO=0
\def\Section#1{\SUBSECNO=0\SUBSUBSECNO=0 \advance\SECNO by 1
        \goodbreak\vskip14pt\noindent{\bf\the\SECNO .\ #1}
        \gdef\Label##1{\xdef##1{\nobreak\the\SECNO}}
        \nobreak\vskip8pt\noindent\kern0pt}
\def\SubSection#1{\SUBSUBSECNO=0 \advance\SUBSECNO by 1
        \goodbreak\vskip14pt\noindent{\it\the\SECNO.\the\SUBSECNO\ #1}
        \gdef\Label##1{\xdef##1{\nobreak\the\SECNO.\the\SUBSECNO}}
        \nobreak\vskip8pt\noindent\kern0pt}
\def\SubSubSection#1{\advance\SUBSUBSECNO by 1
        \goodbreak\vskip14pt\noindent{\rm\the\SECNO.\the\SUBSECNO.\the\SUBSUBSECNO\ #1}
        \gdef\Label##1{\xdef##1{\nobreak\the\SECNO.\the\SUBSECNO.\the\SUBSUBSECNO}}
        \nobreak\vskip8pt\noindent\kern0pt}

\message{definitions;}

\long\def\Definition#1#2{\medbreak\noindent{\bf Definition%
        #1.\enspace}{\it#2}\medbreak\smallskip\relax}
\long\def\Theorem#1#2{\medbreak\noindent{\bf Theorem%
        #1.\enspace}{\it#2}\medbreak\smallskip\relax}
\long\def\Lemma#1#2{\medbreak\noindent{\bf Lemma%
        #1.\enspace}{\it#2}\medbreak\smallskip\relax}

\long\def\Corollary#1#2{\medbreak\noindent{\bf Corollary%
        #1.\enspace}{\it#2}\medbreak\smallskip\relax}

\long\def\Proof#1\Endproof{\noindent{\it Proof.}\enspace
	#1\hfill$\triangleleft$\medbreak\smallskip\relax}

\def\Verbatim{\bgroup%
        \let\Par=\par\def\par{\Par\leavevmode\relax}%
        \obeylines\obeyspaces\tt%
        \lineskip=0pt\baselineskip=0pt%
        \parindent=10pt\def\\{\tt\char092}%
        \catcode`\%=12 \catcode`\#=12 \catcode`\&=12
        \catcode`\$=12 \catcode`\_=12 \catcode`\^=12
        \catcode`\~=12 \catcode`\{=12 \catcode`\}=12
        \relax}
\def\EndVerbatim{\let\par=\Par\vskip-2\parskip\egroup\relax}


\newcount\FOOTNO \FOOTNO=0
\long\def\Footnote#1{\global\advance\FOOTNO by 1
        {\parindent=20pt\baselineskip=9pt%
        \footnote{\nobreak${}^{\the\FOOTNO)}$}{\petit#1\par\vskip-9pt}%
        }\gdef\Label##1{\xdef##1{\nobreak\the\FOOTNO}}}

\message{referencing system: equation numbers,}

\newcount\EQNO \EQNO=0
\def\Eqno{\global\advance\EQNO by 1 \eqno(\the\EQNO)%
        \gdef\Label##1{\xdef##1{\nobreak(\the\EQNO)}}}

\message{figures and tables and}

\def\LArrow[#1]_#2^#3{\mathrel{\smash{\mathop{\hbox to#1{\leftarrowfill}}%
        \limits_{#2}^{#3}}}\relax}
\def\RArrow[#1]_#2^#3{\mathrel{\smash{\mathop{\hbox to#1{\rightarrowfill}}%
        \limits_{#2}^{#3}}}\relax}
\def\UArrow[#1]#2#3{\normallineskip=8pt\matrix{\llap{$\scriptstyle#2$}%
        &\kern-12pt\left\uparrow\vcenter to #1{}\right.&%
        \kern-12pt\rlap{$\scriptstyle#3$}\cr}\relax}
\def\DArrow[#1]#2#3{\normallineskip=8pt\matrix{\llap{$\scriptstyle#2$}%
        &\kern-12pt\left\downarrow\vcenter to #1{}\right.&%
        \kern-12pt\rlap{$\scriptstyle#3$}\cr}\relax}

\newcount\FIGNO \FIGNO=0
\def\Fcaption#1{\global\advance\FIGNO by 1
        {\petit{\bf Fig. \the\FIGNO.~}#1}
        \gdef\Label##1{\xdef##1{\nobreak\the\FIGNO}}}
\def\Figure#1#2{\medskip\vbox{\centerline{\framebox{#1}}
        \centerline{\Fcaption{#2}}}\medskip\relax}

\newcount\TABNO \TABNO=0
\def\Tcaption#1{\global\advance\TABNO by 1
   {\petit{\bf Table. \the\TABNO.~}#1}
   \gdef\Label##1{\xdef##1{\nobreak\the\TABNO}}}

\newread\epsffilein    
\newif\ifepsffileok    
\newif\ifepsfbbfound   
\newif\ifepsfverbose   
\newif\ifepsfdraft     
\newdimen\epsfxsize    
\newdimen\epsfysize    
\newdimen\epsftsize    
\newdimen\epsfrsize    
\newdimen\epsftmp      
\newdimen\pspoints     
\pspoints=1bp          
\epsfxsize=0pt         
\epsfysize=0pt         
\def\epsfbox#1{\global\def\epsfllx{72}\global\def\epsflly{72}%
   \global\def\epsfurx{540}\global\def\epsfury{720}%
   \def\lbracket{[}\def\testit{#1}\ifx\testit\lbracket
   \let\next=\epsfgetlitbb\else\let\next=\epsfnormal\fi\next{#1}}%
\def\epsfgetlitbb#1#2 #3 #4 #5]#6{\epsfgrab #2 #3 #4 #5 .\\%
   \epsfsetgraph{#6}}%
\def\epsfnormal#1{\epsfgetbb{#1}\epsfsetgraph{#1}}%
\def\epsfgetbb#1{%
%
%
\openin\epsffilein=#1
\ifeof\epsffilein\errmessage{I couldn't open #1, will ignore it}\else
%
%
   {\epsffileoktrue \chardef\other=12
    \def\do##1{\catcode`##1=\other}\dospecials \catcode`\ =10
    \loop
       \read\epsffilein to \epsffileline
       \ifeof\epsffilein\epsffileokfalse\else
%
%
          \expandafter\epsfaux\epsffileline:. \\%
       \fi
   \ifepsffileok\repeat
   \ifepsfbbfound\else
    \ifepsfverbose\message{No bounding box comment in #1; using defaults}\fi\fi
   }\closein\epsffilein\fi}%
%
%
\def\epsfclipon{\def\epsfclipstring{ clip}}%
\def\epsfclipoff{\def\epsfclipstring{\ifepsfdraft\space clip\fi}}%
\epsfclipoff
\def\epsfsetgraph#1{%
   \epsfrsize=\epsfury\pspoints
   \advance\epsfrsize by-\epsflly\pspoints
   \epsftsize=\epsfurx\pspoints
   \advance\epsftsize by-\epsfllx\pspoints
%
%
   \epsfxsize\epsfsize\epsftsize\epsfrsize
   \ifnum\epsfxsize=0 \ifnum\epsfysize=0
      \epsfxsize=\epsftsize \epsfysize=\epsfrsize
      \epsfrsize=0pt
%
%
     \else\epsftmp=\epsftsize \divide\epsftmp\epsfrsize
       \epsfxsize=\epsfysize \multiply\epsfxsize\epsftmp
       \multiply\epsftmp\epsfrsize \advance\epsftsize-\epsftmp
       \epsftmp=\epsfysize
       \loop \advance\epsftsize\epsftsize \divide\epsftmp 2
       \ifnum\epsftmp>0
          \ifnum\epsftsize<\epsfrsize\else
             \advance\epsftsize-\epsfrsize \advance\epsfxsize\epsftmp \fi
       \repeat
       \epsfrsize=0pt
     \fi
   \else \ifnum\epsfysize=0
     \epsftmp=\epsfrsize \divide\epsftmp\epsftsize
     \epsfysize=\epsfxsize \multiply\epsfysize\epsftmp
     \multiply\epsftmp\epsftsize \advance\epsfrsize-\epsftmp
     \epsftmp=\epsfxsize
     \loop \advance\epsfrsize\epsfrsize \divide\epsftmp 2
     \ifnum\epsftmp>0
        \ifnum\epsfrsize<\epsftsize\else
           \advance\epsfrsize-\epsftsize \advance\epsfysize\epsftmp \fi
     \repeat
     \epsfrsize=0pt
    \else
     \epsfrsize=\epsfysize
    \fi
   \fi
%
%
   \ifepsfverbose\message{#1: width=\the\epsfxsize, height=\the\epsfysize}\fi
   \epsftmp=10\epsfxsize \divide\epsftmp\pspoints
   \vbox to\epsfysize{\vfil\hbox to\epsfxsize{%
      \ifnum\epsfrsize=0\relax
        \includegraphics{\ifepsfdraft}%
      \else
        \epsfrsize=10\epsfysize \divide\epsfrsize\pspoints
        \includegraphics{\ifepsfdraft}%
      \fi
      \hfil}}%
\global\epsfxsize=0pt\global\epsfysize=0pt}%
%
%
{\catcode`\%=12 \global\let\epsfpercent=
%
%
\long\def\epsfaux#1#2:#3\\{\ifx#1\epsfpercent
   \def\testit{#2}\ifx\testit\epsfbblit
      \epsfgrab #3 . . . \\%
      \epsffileokfalse
      \global\epsfbbfoundtrue
   \fi\else\ifx#1\par\else\epsffileokfalse\fi\fi}%
%
%
\def\epsfempty{}%
\def\epsfgrab #1 #2 #3 #4 #5\\{%
\global\def\epsfllx{#1}\ifx\epsfllx\epsfempty
      \epsfgrab #2 #3 #4 #5 .\\\else
   \global\def\epsflly{#2}%
   \global\def\epsfurx{#3}\global\def\epsfury{#4}\fi}%
%
%
\def\epsfsize#1#2{\epsfxsize}
%
%


\message{the bibliography;}

\newcount\REFNO \REFNO=0
\newbox\REFBOX \setbox\REFBOX=\vbox{}
\def\BegRefs{\message{reading the references}\setbox\REFBOX\vbox\bgroup
        \parindent18pt\baselineskip9pt\petit}
\def\Ref#1{\message{.}\global\advance\REFNO by 1 \ifnum\REFNO>1\vskip3pt\fi
        \item{\the\REFNO .~}\xdef#1{\nobreak[\the\REFNO]}}
\def\EndRefs{\par\egroup\message{done}}
\def\References{\goodbreak\vskip21pt\leftline{\bf References}
        \nobreak\vskip12pt\unvbox\REFBOX\vskip21pt\relax}


\newdimen\defaultparindent
\def\ItemA[#1]#2{\defaultparindent=\parindent\parindent=#1%
        \item{#2}\parindent=\defaultparindent\ignorespaces}%
\def\ItemB#1{\defaultparindent=\parindent\parindent=20pt%
        \item{#1}\parindent=\defaultparindent\ignorespaces}%
\def\ItemOPT{\ifx\ITEMX[\let\ItemCMD\ItemA\else\let\ItemCMD\ItemB\fi\ItemCMD}%
\def\Item{\futurelet\ITEMX\ItemOPT}%

\message{and some math things.  }

\def\N{I\kern-.8ex N}
\def\Z{\raise.72ex\hbox{${}_{\not}$}\kern-.45ex {\rm Z}}
\def\Q{\raise.82ex\hbox{${}_/$}\kern-1.35ex Q} \def\R{I\kern-.8ex R}
\def\C{\raise.87ex\hbox{${}_/$}\kern-1.35ex C} \def\H{I\kern-.8ex H}

\def\RP{{\R\!P}} \def\CP{{\C\!P}} \def\HP{{\H\!P}}

\def\Re{{\rm Re}\,} \def\Im{{\rm Im}\,}


\Date{23.06.1998}
\Title{Transformations of discrete isothermic nets and
	discrete cmc-1 surfaces in hyperbolic space}
	{Udo Hertrich-Jeromin}
	{Dept.\ of Mathematics, Technical University Berlin,
	Str.\ d.\ 17.\ Juni 136, D-10623 Berlin}
	{Using a quaternionic calculus, the Christoffel, Darboux, Goursat,
	and spectral transformations for discrete isothermic nets are
	described, with their interrelations. The Darboux and spectral
	transformations are used to define discrete analogs for cmc-1
	surfaces in hyperbolic space and to obtain a discrete version
	of Bryant's Weierstrass type representation.}
\ShortTitle{Transformations of discrete isothermic nets}
	{U.~Hertrich-Jeromin}


\BegRefs
\Ref\BianchiI L.~Bianchi: {\it Ricerce sulle superficie isoterme
	e sulla deformazione delle quadriche\/}; Annali Mat.~III,
	{\bf 11} (1905) 93-157
\Ref\Bianchi L.~Bianchi: {\it Complementi alle Ricerche sulle
       Superficie Isoterme\/}; Ann.\ Mat.\ pura appl.\ {\bf 12}
       (1905) 19--54
\Ref\Blaschke W.~Blaschke: {\it Vorlesungen \"uber Differentialgeometrie
        III\/}; Springer, Berlin 1929
\Ref\BoJe A.~Bobenko, U.~Hertrich-Jeromin: {\it Orthogonal nets and
	Clifford algebras\/}; Preprint, 1998
\Ref\BoPiA A.~Bobenko, U.~Pinkall: {\it Discrete Isothermic Surfaces\/};
        J.~reine angew.\ Math. {\bf 475} (1996) 187-208
\Ref\BoPiB A.~Bobenko, U.~Pinkall: {\it Discretization of Surfaces and
        Integrable systems\/}; in A.~Bobenko, R.~Seiler: Discrete
        Integrable Geometry and Physics, Oxford UP, Oxford 1999
\Ref\BJPP F.~Burstall, U.~Hertrich-Jeromin, F.~Pedit, U.~Pinkall:
	{\it Curved flats and Isothermic surfaces}; Math.~Z.\ {\bf 225}
	(1997) 199-209
\Ref\Fran F.~Burstall: {\it Isothermic surfaces: Conformal geometry,
	Clifford algebras, and Integrable systems\/}; in preparation
\Ref\BPP F.~Burstall, F.~Pedit, U.~Pinkall: {\it Isothermic submanifolds
	of a symmetric $R$-space\/}; in preparation
\Ref\Cartan \'{E}.~Cartan: {\it Les espaces \`{a} connexion conforme\/};
        Ann.\ Soc.\ Pol.\ Math.\ {\bf 2} (1923) 171-221
\Ref\Christoffel E.~Christoffel: {\it Ueber einige allgemeine
        Eigenschaften der Minimumsfl\"achen\/}; Crelle's J.\ {\bf 67}
        (1867) 218--228
\Ref\Darboux G.~Darboux: {\it Sur les surfaces isothermiques\/};
        Comptes Rendus {\bf 122} (1899) 1299--1305, 1483--1487, 1538
\Ref\Suppl U.~Hertrich-Jeromin: {\it Supplement on Curved flats in the space
        of Point pairs and Isothermic surfaces: A Quaternionic calculus\/};
        Doc.\ Math.\ J.~DMV {\bf 2} (1997) 335-350
\Ref\JePe U.~Hertrich-Jeromin, F.~Pedit: {\it Remarks on the Darboux
        transform of isothermic surfaces\/};
        Doc.\ Math.\ J.~DMV {\bf 2} (1997) 313-333
\Ref\JMN U.~Hertrich-Jeromin, E.~Musso, L.~Nicolodi: {\it M\"obius geometry
	of surfaces of constant mean curvature 1 in hyperbolic space\/};
	Preprint, 1998
\Ref\HHP U.~Hertrich-Jeromin, T.~Hoffmann, U.~Pinkall: {\it A discrete
        version of the Darboux transform for isothermic surfaces\/};
        in A.~Bobenko, R.~Seiler, {\it Discrete integrable Geometry and
	Physics\/}, Oxford Univ.\ Press, Oxford 1999
\Ref\Musso E.~Musso: {\it Deformazione di superficie nello spazio
        di M\"obius\/}; Rend.\ Istit.\ Mat.\ Univ.\ Trieste {\bf 27}
        (1995) 25-45
\Ref\PePi F.~Pedit, U.~Pinkall: {\it Rigidity of surfaces in M\"obius
	geometry\/}; personal communication
\Ref\Sauer R.~Sauer: {\it Differenzengeometrie\/}; Springer, Berlin 1970
\Ref\UmYa M.~Umehara, K.~Yamada: {\it A parametrization for the Weierstrass
	formulae and perturbation of complete minimal surfaces in $\R^3$ into
	the hyperbolic 3-space\/}; J.~reine angew.\ Math.\ {\bf 432} (1997)
	93-116
\Ref\UmYaD M.~Umehara, K.~Yamada: {\it A Duality on cmc-1 surfaces in
	hyperbolic space and a hyperbolic analogue of the Ossermann
	inequality\/}; Tsukuba J.~Math.\ {\bf 21} (1997) 229-237
\EndRefs

\def\CatenoidCousins{\Figure{\hbox{
        \epsfxsize=.32\hsize\epsfclipon\epsfbox[110 230 490 610]{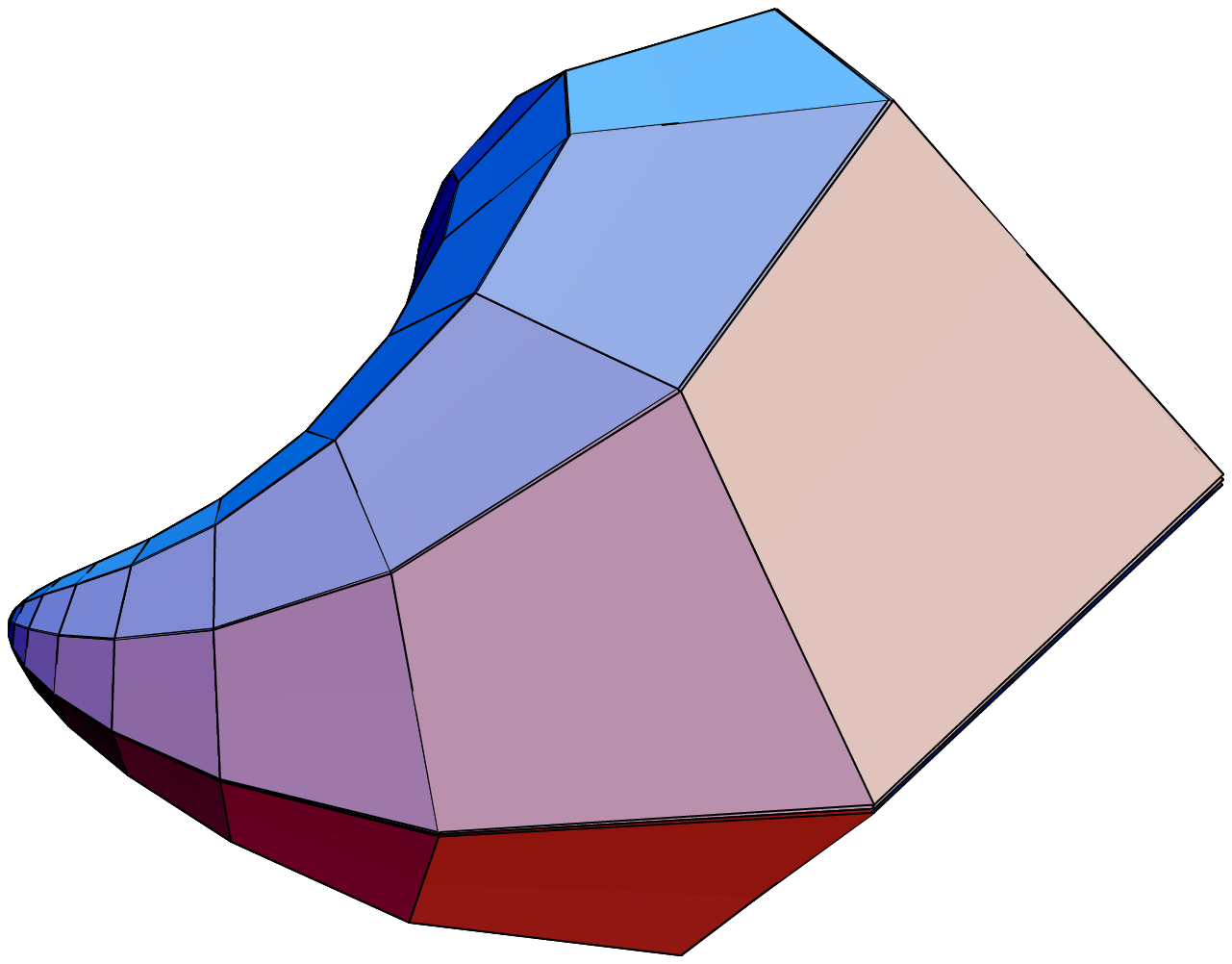}
        \vrule
        \epsfxsize=.32\hsize\epsfclipon\epsfbox[125 295 455 625]{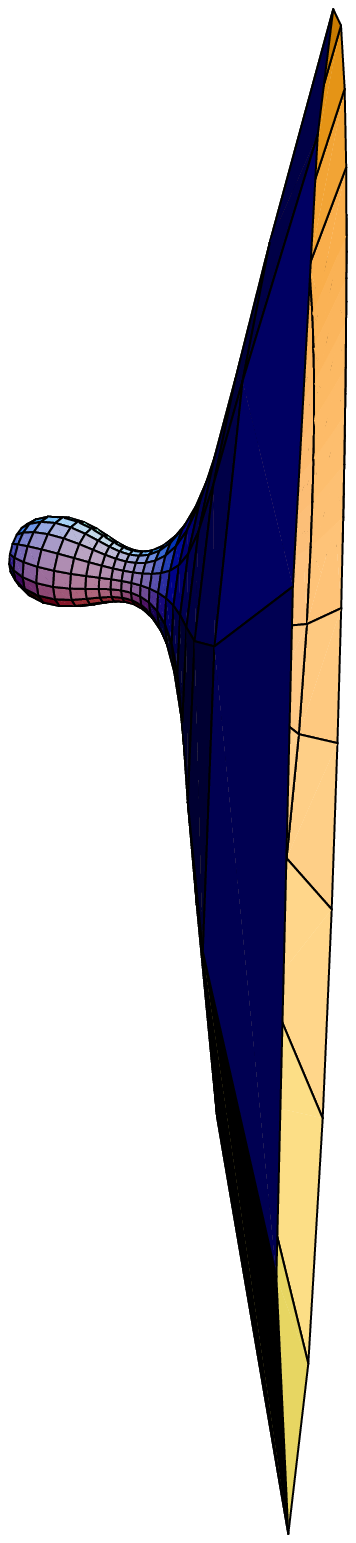}
        \vrule
        \epsfxsize=.32\hsize\epsfclipon\epsfbox[120 190 535 615]{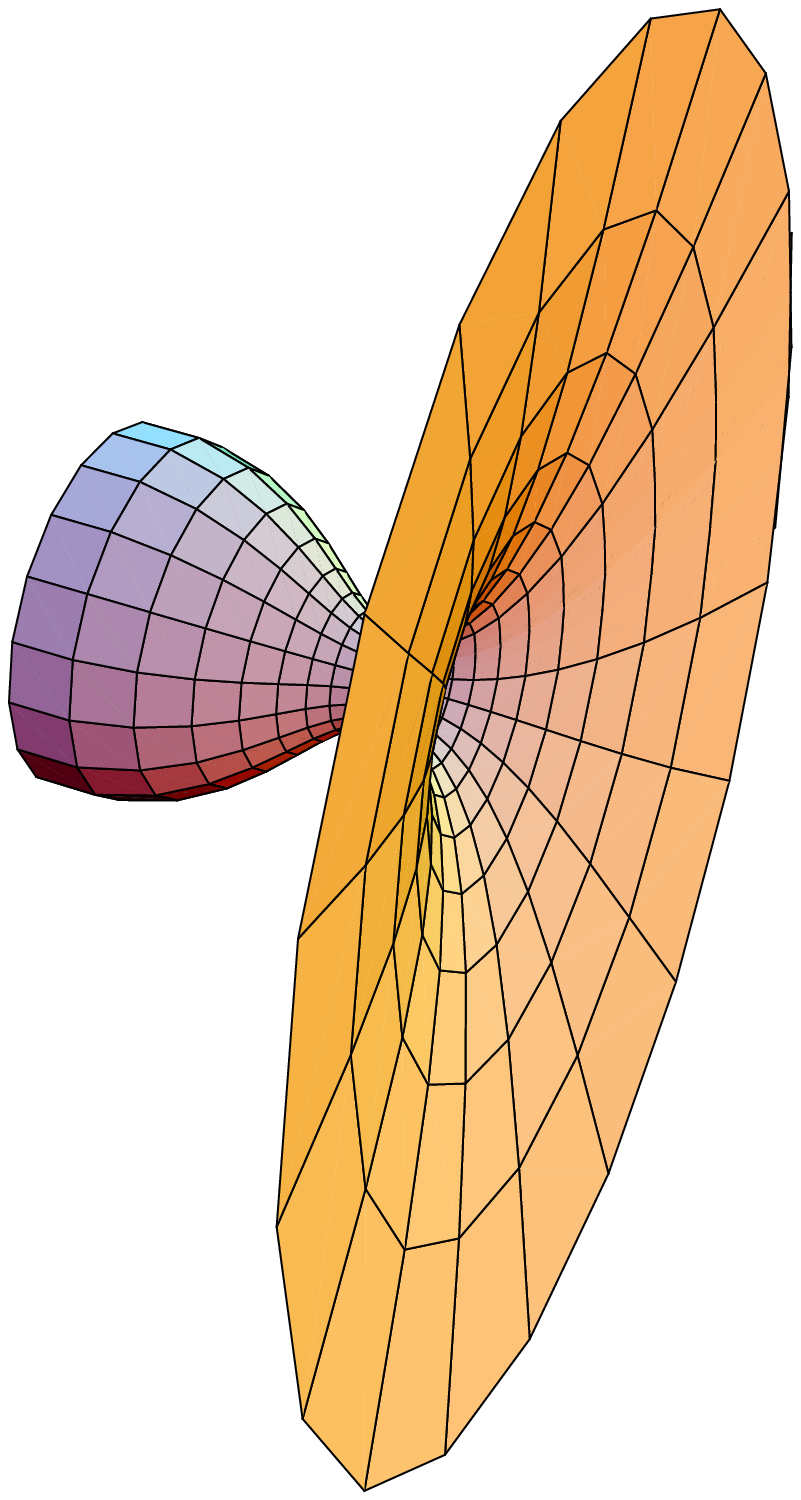}
        }\hrule\hbox{
        \epsfxsize=.32\hsize\epsfclipon\epsfbox[115 170 550 615]{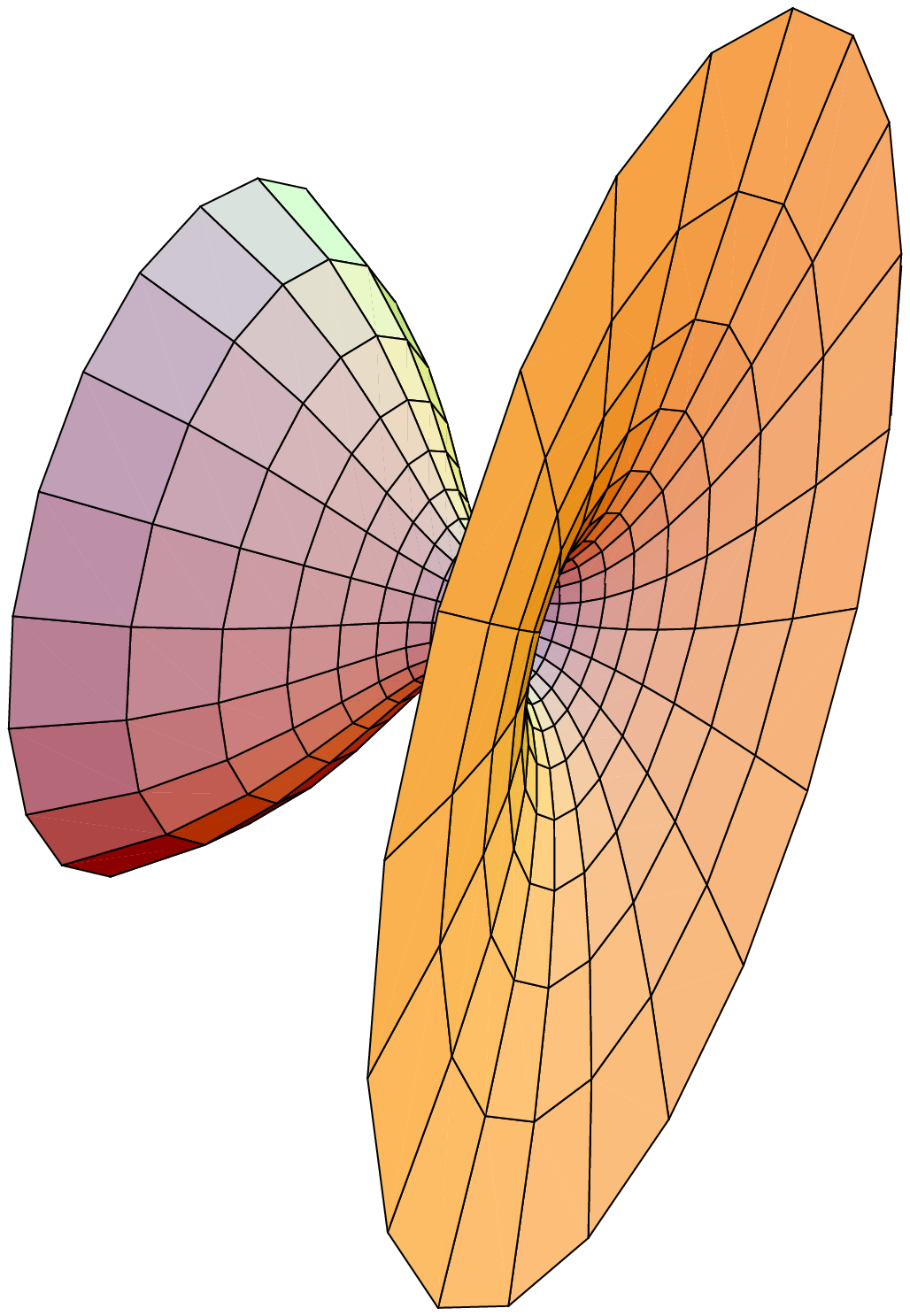}
        \vrule
        \epsfxsize=.32\hsize\epsfclipon\epsfbox[072 170 547 645]{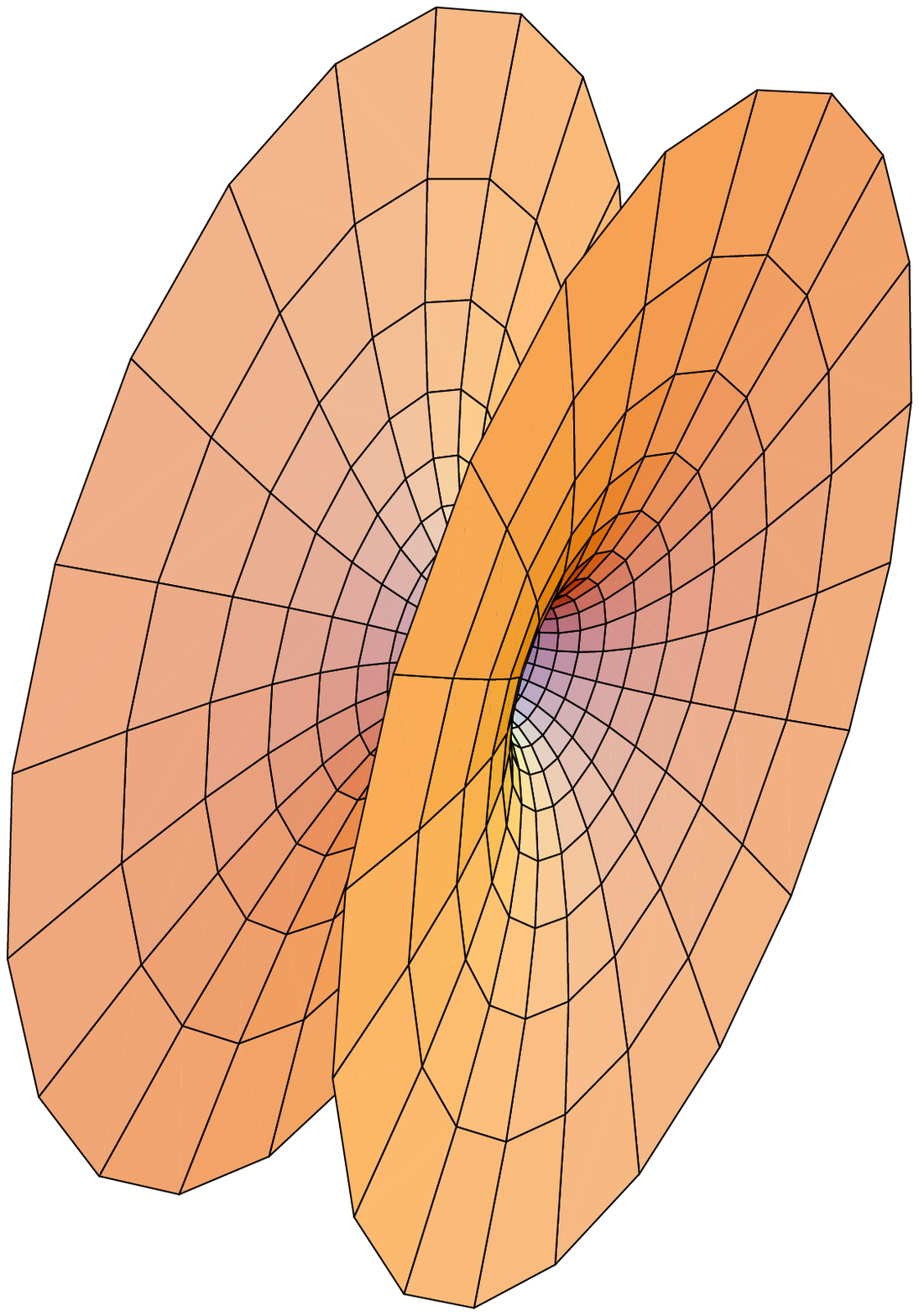}
        \vrule
        \epsfxsize=.32\hsize\epsfclipon\epsfbox[095 180 545 630]{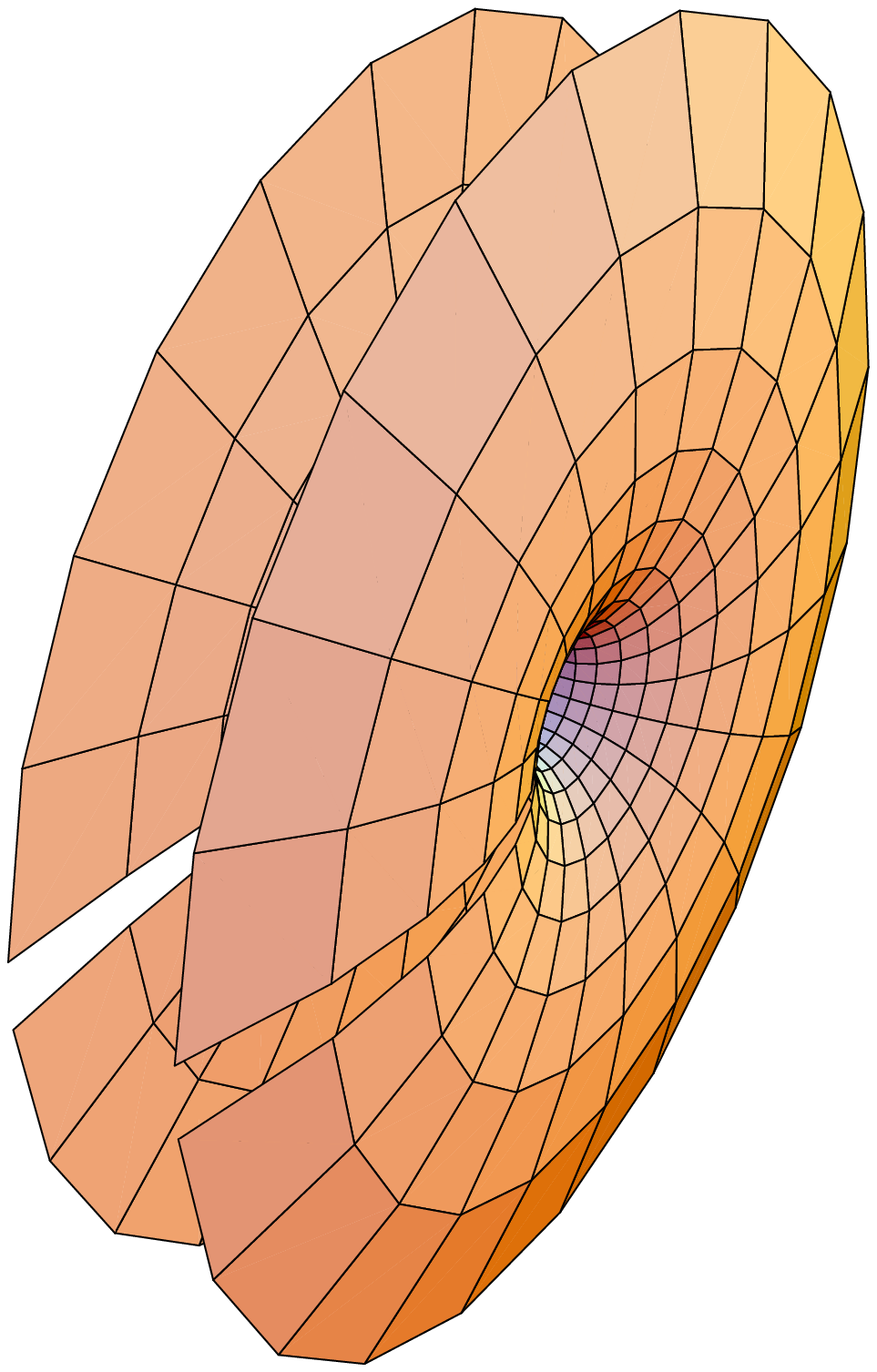}
        }\hrule\hbox{
        \epsfxsize=.32\hsize\epsfclipon\epsfbox[100 200 540 640]{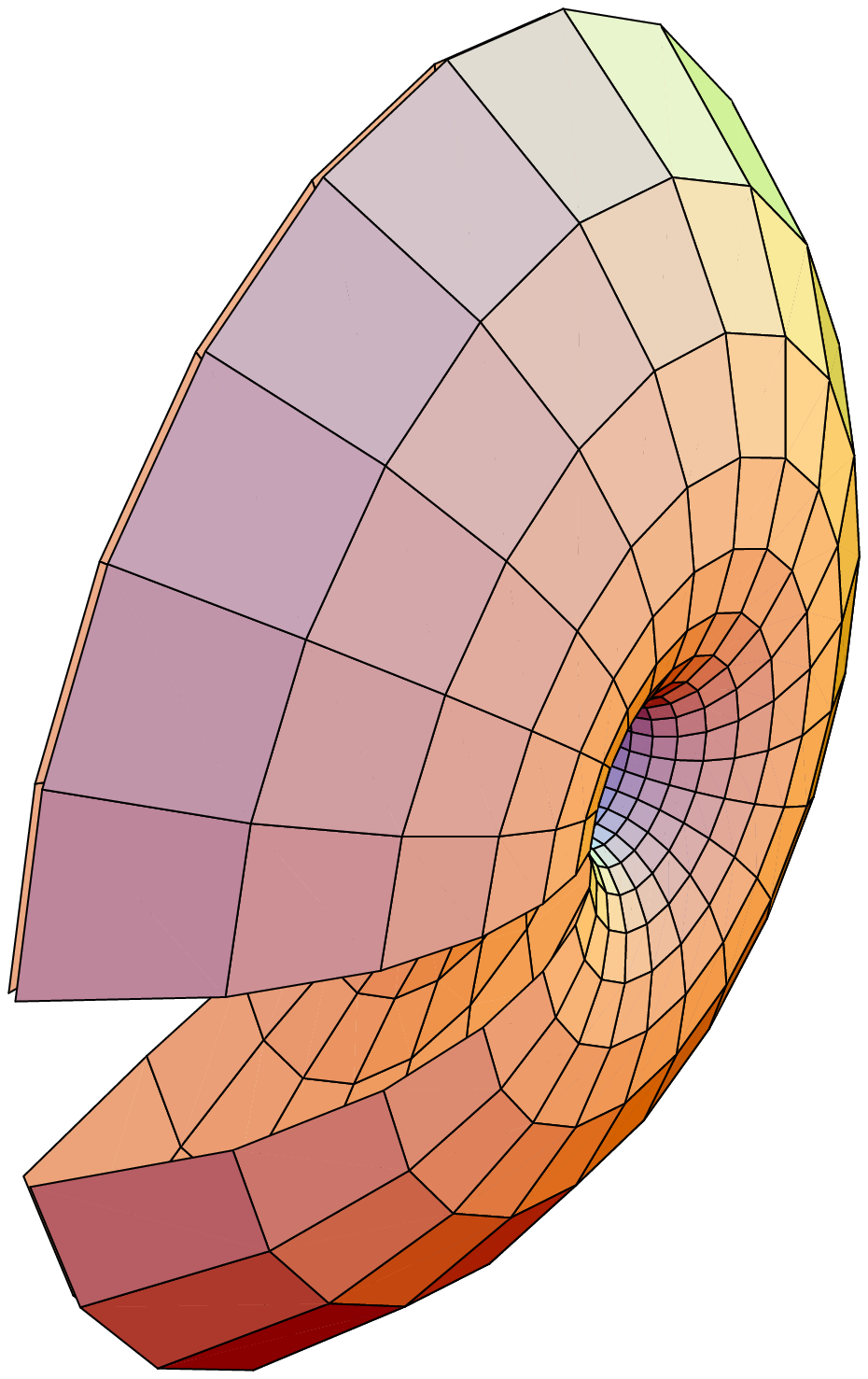}
        \vrule
        \epsfxsize=.32\hsize\epsfclipon\epsfbox[110 195 545 630]{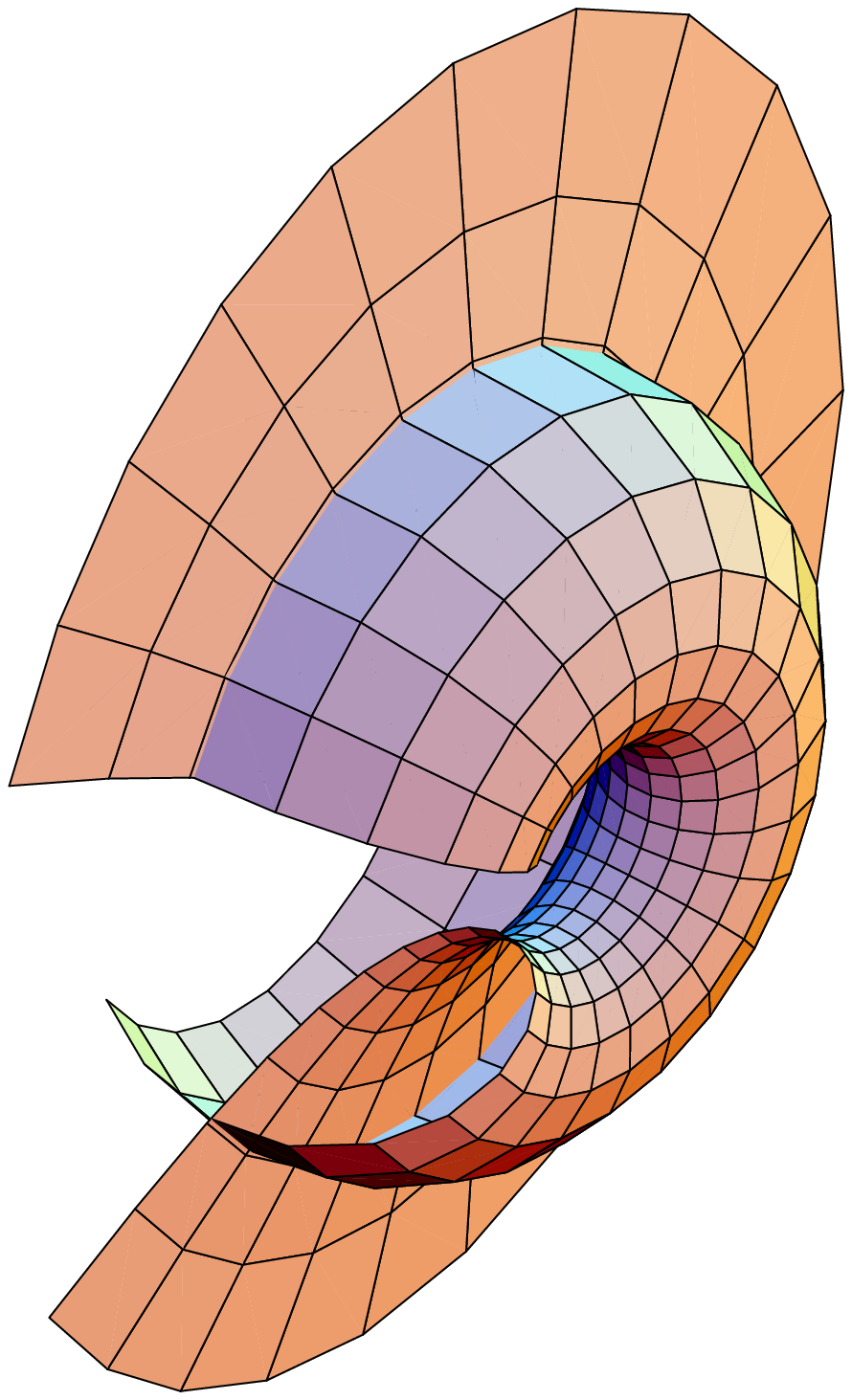}
        \vrule
        \epsfxsize=.32\hsize\epsfclipon\epsfbox[075 195 505 625]{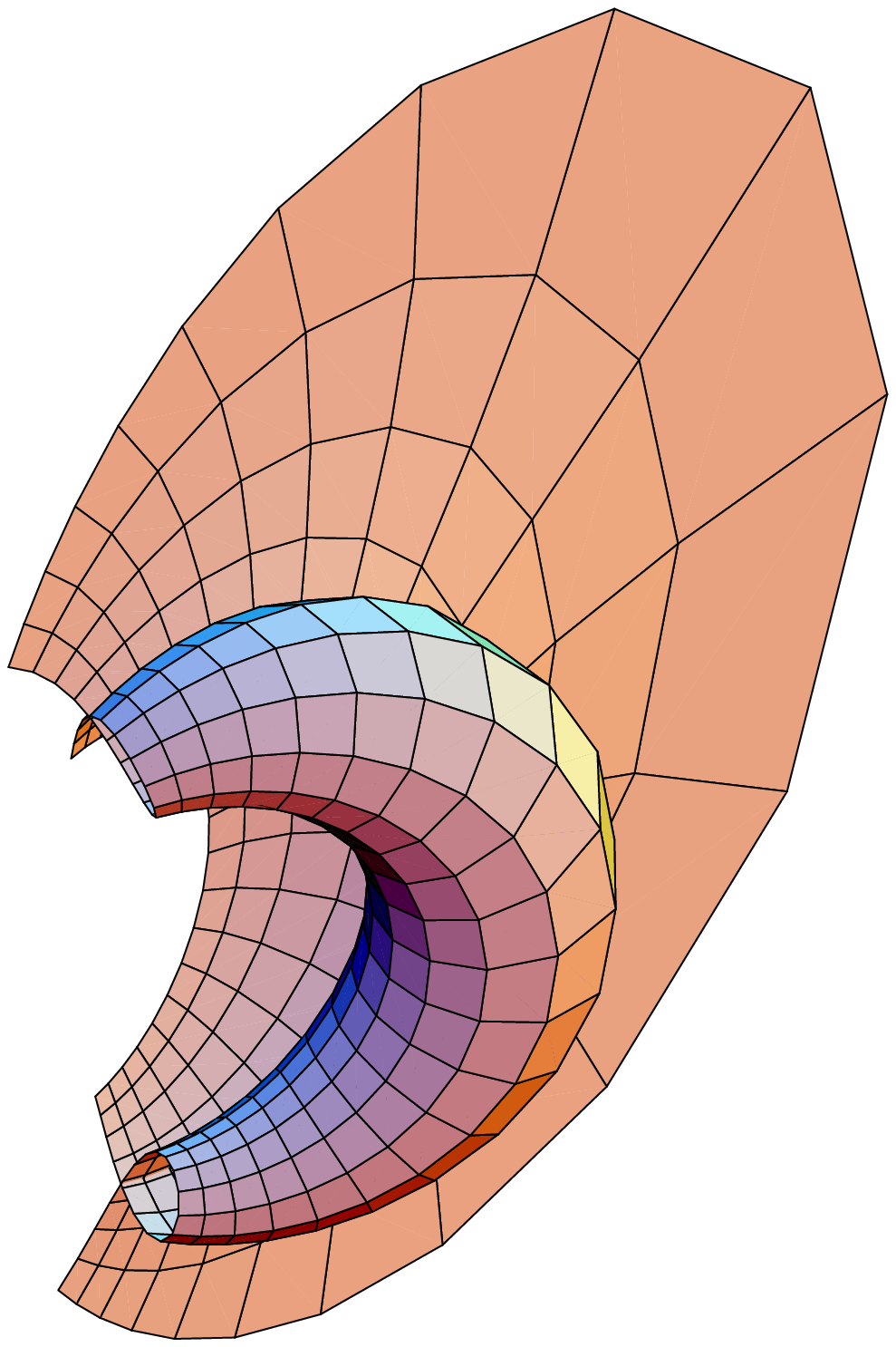}
        }
        }{$T$-transforms of a discrete Catenoid with Gauss map
        ${\fk g}_{m,n}=e^{2\pi{m+in\over20}}$}\Label\CatCous}

\Section{Introduction}
There is a rich transformation theory for smooth isothermic surfaces
that was mostly developed in classical times:\par
	Christoffel (1867) showed that, to a given isothermic surface
(in Euclidean 3-space), there is (locally) a second isothermic surface
and a conformal point-to-point correspondence between the two surfaces
such that the tangent planes in corresponding points are parallel ---
in fact, he showed that the existence of such a (non-trivial) transform
{\it characterizes} isothermic surfaces \Christoffel;\par
	Darboux (1899) found the existence of a $(1+3)$-parameter family
of sphere congruences that are enveloped by a given isothermic surface,
such that the induced map onto the second envelope preserves curvature
lines and is conformal --- this second envelope (transform) is then also
isothermic \Darboux;\par
	Bianchi (1905) introduced a 1-parameter family of deformations
for an isothermic surface --- that he called the ``$T$-transformation''
--- and analyzed the interrelations of the various transformations in
terms of ``permutability theorems'' \Bianchi.
This last transformation seems to play a central role in the theory of
isothermic surfaces: it relates to the characterization of isothermic
surfaces as the only non-rigid surfaces in M\"obius geometry \Musso\
(cf.\Cartan,\PePi) as well as to the integrable system approach to
isothermic surfaces via curved flats in the symmetric space of point
pairs \BJPP\ (cf.\Suppl,\BPP).

Discretizations for the Christoffel- and the Darboux-transformations
have been given, using a discrete version of Christoffel's original
formula and of the Riccati-type partial differential equation that we
derived in \JePe: for both ansatzes, it was crucial to have a Euclidean
description of the transformation (even though, for the discrete analog
of the Darboux-transformation, a M\"obius invariant definition can be
given --- as it should be).
The $T$-transformation is rather a transformation for M\"obius equivalence
classes of isothermic surfaces than for surfaces with a fixed position in
space (in contrast to the Darboux transformation that incorporates the
position of the surface); consistently, any desription should make use
of the M\"obius group. Thus, in order to introduce a discrete analog for
the $T$-transformation, it seemed necessary to elaborate a suitable setup
(section 2): a quaternionic approach\Footnote{Or, in higher dimensions, an
approach using Clifford algebras \Fran\ (or, \BoJe).} seems preferable over
the classical formalism as no normalization problems occur.

With this setup, the transformation theory of discrete isothermic nets is
elaborated (section 3) --- entirely analogous to the transformations of
(smooth) isothermic surfaces (cf.\JMN): besides the $T$-transformation,
a notion of ``general Christoffel ($C$-) transformation'' (cf.\BPP) and
a (generalized) ``Goursat ($G$-) transformation'' are introduced; also, a
new definition for the Darboux ($D$-) transformation is given and shown to
be equivalent to our previous definition \HHP. Finally, the interrelations
between the different transformations are discussed in terms of permutability
theorems.

The transformations of isothermic surfaces already provided useful ansatzes
to define discrete minimal (cf.\BoPiA) and discrete cmc \HHP\ nets.
With the $T$-transformation available, we may use it to define discrete analogs 
for surfaces of constant mean curvature in other space forms: in case of cmc
surfaces, the $T$-transformation becomes Lawson correspondence (cf.\Suppl)
--- and, in particular, it becomes the Umehara-Yamada perturbation in case
of minimal surfaces in Euclidean space and cmc-1 surfaces in hyperbolic space
(cf.\JMN). Thus, discrete analogs for surfaces of constant mean curvature $H$
in a space of constant curvature $k$, that can be obtained from cmc surfaces in
Euclidean space via Lawson correspondence (i.e.\ $k+H^2\geq0$), can be defined
as $T$-transforms of minimal or cmc nets in Euclidean space.

However, for cmc-1 surfaces in hyperbolic space, a more geometric approach is
obtained via the Darboux transformation as it couples the surface with its
hyperbolic Gauss map \JMN: this approach will be used to introduce the notion
of ``horospherical nets'' as an analog of cmc-1 surfaces in hyperbolic space,
in section 4. The permutability theorems for various transformations will
show that this ansatz is equivalent to the previously described one --- this
way, the notion of the ``minimal cousin'' of a horospherical net makes sense.
Also, they will serve in obtaining a discrete version of Bryant's Weierstrass
type representation;
comparison with the other Weierstrass type representation in terms of the
Darboux transformation, a discrete analog of the ``dual transformation''
for cmc-1 surfaces by Umehara and Yamada \UmYaD\ naturally appears.

\Section{Preliminaries}
Just as the geometry of the complex projective line $\CP^1\cong S^2$
is the geometry of (orientation preserving) M\"obius transformations
of the conformal 2-sphere, 1-dimensional quaternionic projective
geometry can be identified with the M\"obius geometry of the conformal
4-sphere $S^4\cong\HP^1:=\{v\H\,|\,v\in\H^2\}$ (cf.\Suppl).
Here, we consider the space $\H^2$ of homogeneous coordinates as a
{\it right} vector space over the quaternions.
In this model, $PGl(2,\H)=Gl(2,\H)/\R I$ acts by (orientation preserving)
M\"obius transformations on $\HP^1$ via $v\H\to(Mv)\H$, and
3-spheres appear as null cones of suitable quaternionic hermitian
forms $s:\H^2\times\H^2\to\H$: equipped with the Minkowski product
induced\Footnote{We set $\det s:=s(e_1,e_1)s(e_2,e_2)-|s(e_1,e_2)|^2$
relative to a fixed basis $(e_1,e_2)$ of $\H^2$ --- a different choice
of basis results in a rescaling of the Minkowski product: this ambiguity
reflects the fact that the geometrically relevant space is the projective
space $\RP^5$ with absolute quadric $\det=0$.}
by $-\det$ as a quadratic form the space $\R^6_1$ of quaternionic
hermitian forms becomes the space of homogeneous coordinates for the
classical model space $\RP^5$ of 4-dimensional M\"obius geometry
(\Blaschke, cf.\Suppl) with $PGl(2,\H)$ acting via $(M,s)\to M\cdot s
:=s(M^{-1}.,M^{-1}.)$.

To any quaternionic line $v\H\subset\H^2$ there is a unique line $\H\nu
=(v\H)^{\perp}\subset(\H^2)^{\ast}$ of quaternionic 1-forms annihilating
$v\H$ --- note that the dual space $(\H^2)^{\ast}$ is a {\it left} vector
space over $\H$ since $\H^2$ is a right vector space.
Thus, the conformal 4-sphere can be identified $(\HP^1)^{\perp}:=\{\H\nu\,|\,
\nu\in(\H^2)^{\ast}\}\cong\HP^1$.  And, $PGl(2,\H)$ acts on $(\HP^1)^{\perp}$
via $\H\nu\to\H(M\cdot \nu)=\H(\nu M^{-1})$.

Fixing homogeneous coordinates $s\in\R^6_1$ for a 3-sphere $S^3\cong\R
s\subset\R^6_1$ 
provides a canonical identification $\H^2\leftrightarrow(\H^2)^{\ast}$
via $v\to s(v,.)$.
Indeed, $v\to\sigma_v:=s(v,.)$ is an anti-isomorphism as it is anti-linear
and $\sigma_v=0$ if and only if $v=0$: $$\matrix{ |s|^2v
	&=& e_1\{s(e_1,e_2)s(e_2,v)-s(e_1,v)s(e_2,e_2)\}\hfill\cr
	&+& e_2\{s(e_2,e_1)s(e_1,v)-s(e_2,v)s(e_1,e_1)\}.\cr}
	\Eqno\Label\KramerRule$$
Moreover, for $v\in S^3$, i.e.\ $s(v,v)=0$ and $v\neq0$, $(v\H)^{\perp}
=\H\sigma_v$ since $\sigma_v\neq0$ implies that $\sigma_v(w)\neq0$ unless
$w=v\lambda$. Consequently, the quaternionic hermitian form $s^{\ast}:
(\H^2)^{\ast}\times(\H^2)^{\ast}\to\H$, $s^{\ast}(\sigma_v,\sigma_w)=s(v,w)$
defines $S^3$ in the $(\HP^1)^{\perp}$-model of $S^4$.
The M\"obius transformations of $S^3\cong\R s$ are those fixing $\R s$, i.e.\
$\R(M^{-1}\cdot s)=\R s$.

The above identification $\HP^1\cong(\HP^1)^{\perp}$ can be used to
calculate the cross ratio \HHP\ of four points in $\HP^1$, or in $S^3$,
respectively:

\Lemma{~(Cross ratio)}{The cross ratio of four points $p_i\in\HP^1$, $i=1,
\dots,4$, is given by $[p_1,p_2,p_3,p_4]=\Re q+i\,|\Im q|$, where $$q=
(\nu_1v_2)(\nu_3v_2)^{-1}(\nu_3v_4)(\nu_1v_4)^{-1}\Eqno\Label\CrossRatio$$
with homogeneous coordinates $v_i\in\H^2$ and $\nu_i\in(\H^2)^{\ast}$ of
$p_i\in\HP^1\cong(\HP^1)^{\perp}$.\par
For four points $p_i\in S^3\cong s$ their homogeneous coordinates $\nu_i=
\sigma_{v_i}\simeq v_i$ in $(\H^2)^{\ast}\cong\H^2$ can be identified, so
that $$q=(s(v_1,v_2))(s(v_3,v_2))^{-1}(s(v_3,v_4))(s(v_1,v_4))^{-1}.\Eqno$$}

Obviously, \CrossRatio\ is independent of the choice of homogeneous
coordinates, as well as it is invariant under M\"obius transformations.
The easiest way to see that \CrossRatio\ really gives the cross ratio
of the four points is to use affine coordinates: without loss of generality
$p_i\neq\pmatrix{1\cr0\cr}\H\simeq\infty\in\HP^1$ so that one can choose
$v_i=\pmatrix{{\fk p}_i\cr1\cr}$ and $\nu_i=(1,-{\fk p}_i)$ with suitable
${\fk p}_i\in\H$; this way, one obtains the usual (cf.\HHP) formula $q=
({\fk p}_1-{\fk p}_2)({\fk p}_2-{\fk p}_3)^{-1}({\fk p}_3-{\fk p}_4)
({\fk p}_4-{\fk p}_1)^{-1}$.

A more invariant way of describing the choice of affine coordinates is
via stereo\-graphic projection:
choosing homogeneous coordinates $v_0,v_{\infty}\in\H^2$ of two points
$0,\infty\in\HP^1\cong(\HP^1)^{\perp}$ their homogeneous coordinates
$\nu_0,\nu_{\infty}\in(\H^2)^{\ast}$ can be canonically fixed by the
condition $v_0\nu_{\infty}+v_{\infty}\nu_0=id$, i.e.\ $(v_0,v_{\infty})$
and $(\nu_0,\nu_{\infty})$ are ``pseudo dual'' bases.
Then, if $v\in\H^2$ and $\nu\in(\H^2)^{\ast}$ denote homogeneous coordinates
of a point $\infty\neq p=\HP^1\cong(\HP^1)^{\perp}$, $v\H=(v_{\infty}{\fk p}
+v_0)\H$ and $\H\nu=\H(\nu_0-{\fk p}\nu_{\infty})$ with ${\fk p}:=(\nu_0v)
(\nu_{\infty}v)^{-1}=-(\nu v_{\infty})^{-1}(\nu v_0)$ where the last equality
follows from $\nu v=0$.
${\fk p}={\fk p}(p)$ will be called the ``stereographic projection'' of
$p\in\HP^1$. Note that ${\fk p}$ does not depend on the choice of homogeneous
coordinates $v$, or $\nu$, of $p$;
a different choice of homogeneous coordinates for $0$ and $\infty$, $\tilde
v_0=v_0a_0$ and $\tilde v_{\infty}=v_{\infty}a_{\infty}^{-1}$ or, equivalently,
$\tilde\nu_0=a_{\infty}\nu_0$ and $\tilde\nu_{\infty}=a_0^{-1}\nu_{\infty}$,
yields a stretch-rotation $\tilde{\fk p}=a_{\infty}{\fk p}a_0$ of the
stereographic projection.
Now, for two points $p_1,p_2\in\HP^1\cong(\HP^1)^{\perp}$, a simple
calculation yields $$(\nu_1v_{\infty})^{-1}(\nu_1v_2)(\nu_{\infty}v_2)^{-1}
={\fk p}_2-{\fk p}_1.\Eqno\Label\Edge$$

Restricting to points in a 3-sphere $S^3\cong\R s$, the identification
$\H^2\cong(\H^2)^{\ast}$ can be used to fix the relative scaling of
homogeneous coordinates $v_0,v_{\infty}\in\H^2$ for $0,\infty\in
S^3\subset\HP^1$, via $s(v_0,v_{\infty})=1$; then, $|s|=1$ with respect
to $(v_0,v_{\infty})$ as a basis of $\H^2$, and \KramerRule\ reduces to $v=
v_0s(v_{\infty},v)+v_{\infty}s(v_0,v)=(v_0+v_{\infty}{\fk p})s(v_{\infty},v)$.
Hence, $s(v,v)=0$ implies ${\fk p}+\bar{\fk p}=0$, such that the stereographic
projection takes values in the imaginary quaternions.
A different choice of homogeneous coordinates $\tilde v_0=v_0a$ now determines
$\tilde v_{\infty}=v_{\infty}\bar a^{-1}$ --- the effect on the stereographic
projection being a stretch-rotation $\tilde{\fk p}=a{\fk p}\bar a$ in
$\R^3\cong{\rm Im}\H$.


\Section{Transformations of discrete isothermic nets}\Label\Conventions
Now, let $f:\Z^2\to\HP^1$ denote a discrete map, $\varphi$ the corresponding
map into $(\HP^1)^{\perp}$, i.e.\ $\varphi\circ f=0$, and let ${\fk f}:\Z^2\to
\H$ denote its stereographic projection with respect to $(v_0,v_{\infty})$, or
its pseudo dual basis $(\nu_0,\nu_{\infty})$ respectively.
To simplify formulas, we will assume $\varphi v_{\infty}=\nu_{\infty}f\equiv1$
wherever a choice of homogeneous coordinates for two points $0,\infty$ plays
role. Then \Edge\ yields $$\matrix{
	-\varphi_{m+1,n}f_{m,n} &=& \varphi_{m,n} f_{m+1,n}
		&=& {\fk f}_{m+1,n} - {\fk f}_{m,n}
		&=:& (\partial_1{\fk f})_{m,n}, \cr
	-\varphi_{m,n+1}f_{m,n} &=& \varphi_{m,n} f_{m,n+1}
		&=& {\fk f}_{m,n+1} - {\fk f}_{m,n}
		&=:& (\partial_2{\fk f})_{m,n}. \cr
	}\Eqno\Label\Edges$$
Thus, a minimal regularity assumption on $f$ is suggested:

\Definition{~(Discrete net)}{A discrete map $\varphi\simeq f:\Z^2\to\HP^1
\cong(\HP^1)^{\perp}$, or ${\fk f}:\Z^2\to\H$ respectively, will be called a
discrete net if $(\partial_1{\fk f})_{m,n}\neq0\neq(\partial_2{\fk f})_{m,n}$
and $(\partial_1{\fk f})_{m,n}\not\,\parallel(\partial_2{\fk f})_{m,n}$ are
non-parallel for all $(m,n)\in\Z^2$.}

For the rest of the paper, we will adopt the above notation conventions;
in case of discrete nets $f:\Z^2\to S^3\cong\R s\subset\R^6_1$ into the
conformal 3-sphere, we will assume $s(v_0,v_{\infty})=1$ and $\varphi=s(f,.)$.

As an analog of an arbitrary (not necessarily conformal) curvature line
para\-metrization of a smooth isothermic surface we use the following
``wide definition'' (cf.\BoPiB) for discrete isothermic nets:

\Definition{~(Isothermic net)}{A discrete net $f:\Z^2\to\HP^1$ is called
isothermic if all cross ratios
$q_{m,n}=[f_{m,n},f_{m+1,n},f_{m+1,n+1},f_{m,n+1}]={a_m\over b_n}$
where the (real valued) functions $a,b:\Z\to\R$ depend on $m$ resp.\ $n$ only.}

This way, isothermic nets are a special case of discrete principal nets \BoPiB:

\Definition{~(Principal net)}{A discrete net $f:\Z^2\to\HP^1$ is called
principal net if all elementary quadrilaterals are concircular, i.e.\ if
$q:\Z^2\to\R$.}

\SubSection{The difference equations}
Smooth isothermic surfaces $f:M^2\to S^3$ can be characterized as the only
``M\"obius deformable'' surfaces, i.e.\ by the existence of a family of
non constant maps $T^{\lambda}:M^2\to M\ddot{o}b(S^3)$ such that $f^{\lambda}
:=T^{\lambda}f:M^2\to S^3$ induce the same curvature lines and conformal
metric\Footnote{The term ``conformal metric'' is used for the induced metric
of the central sphere congruence (conformal Gauss map) of the surface; the
``conformal Hopf differential'' of a surface encodes both informations,
thus it determines an immersion $f$ uniquely up to M\"obius transformation
unless $f$ parametrizes an isothermic surface (cf.\PePi).} (cf.\Cartan,\Musso).
Classically, this deformation is known as the ``$T$-transformation'' \Bianchi;
we first give a discrete version of the corresponding system of differential
equations in the quaternionic setting \BPP,\JMN:

{\Lemma{~($T$-transformation)}{Let $f:\Z^2\to\HP^1$, $f\simeq\varphi:\Z^2\to
(\HP^1)^{\perp}$, be a discrete net, and $\nu_{\infty}f=\varphi v_{\infty}
\equiv1$ after some choice of homogeneous coordinates $v_0,v_{\infty}\in\H^2$
for two points $0,\infty\in\HP^1$; we define $$\matrix{
	U_{m,n} &:=& f_{m,n}u_{m,n}\varphi_{m+1,n} &&{\rm and}&&
	V_{m,n} &:=& f_{m,n}v_{m,n}\varphi_{m,n+1}\cr}\Eqno\Label\RForm$$
with some functions $u,v:\Z^2\to\H$. Then, the system
$$\matrix{
	T^{\lambda}_{m+1,n} &=& T^{\lambda}_{m,n}(1+\lambda U_{m,n}) &,&
	T^{\lambda}_{m,n+1} &=& T^{\lambda}_{m,n}(1+\lambda V_{m,n}) \cr}
	\Eqno\Label\TSystem$$
is solvable for all $\lambda\in\R$ if and only if the following system is:
$$\matrix{
	F^{\lambda}_{m+1,n} &=& F^{\lambda}_{m,n}(1
		+v_{\infty}(\partial_1{\fk f})_{m,n}\nu_{\infty}
		+\lambda v_0u_{m,n}\nu_0), \cr
	F^{\lambda}_{m,n+1} &=& F^{\lambda}_{m,n}(1
		+v_{\infty}(\partial_2{\fk f})_{m,n}\nu_{\infty}
		+\lambda v_0v_{m,n}\nu_0). \cr
	}\Eqno\Label\OldTSystem$$}\PLabel\TLemma
\Proof
Since $\varphi v_{\infty}=\nu_{\infty}f\equiv1$, we have $\nu_0f={\fk f}
=-\varphi v_0$ as $\varphi f\equiv0$.  Thus, with the ansatz
$F^0:=id+v_{\infty}{\fk f}\,\nu_{\infty}$ for the ``Euclidean frame''
of $f$, a straightforward calculation yields
$$\matrix{
	(F^0_{m,n})^{-1}(1+\lambda U_{m,n})\,F^0_{m+1,n} &=& (id
                +v_{\infty}(\partial_1{\fk f})_{m,n}\nu_{\infty}
                +\lambda v_0u_{m,n}\nu_0), \cr
	(F^0_{m,n})^{-1}(1+\lambda V_{m,n})\,F^0_{m,n+1} &=& (id
                +v_{\infty}(\partial_2{\fk f})_{m,n}\nu_{\infty}
                +\lambda v_0v_{m,n}\nu_0), \cr}$$
i.e.\ the gauge equivalence of the systems \TSystem\ and \OldTSystem.
\Endproof}

This lemma gives the relation between the discrete ``curved flat''
system \OldTSystem\ and the system \TSystem\ modelling the conformal
deformability of (isothermic) nets:
using \Edges\ and denoting $a=u\cdot(\partial_1{\fk f})$ and
$b=v\cdot(\partial_2{\fk f})$, we have $$\matrix{
	(1+\lambda U_{m,n})f_{m,n} \hfill&=&
		f_{m,n}(1-\lambda a_{m,n}) , \hfill\cr
	(1+\lambda U_{m,n})f_{m+1,n} \hfill&=& f_{m+1,n}; \hfill\cr
	(1+\lambda V_{m,n})f_{m,n} \hfill&=&
		f_{m,n}(1-\lambda b_{m,n}), \cr
	(1+\lambda V_{m,n})f_{m,n+1} \hfill&=& f_{m,n+1}. \hfill\cr }
	\Eqno\Label\TFixpoints$$
Thus, the $T^{\lambda}$ (if they exist) map a ``vertex star'' consisting
of a center vertex $f_{m,n}$ and its four neighbours $f_{m\pm1,n}$ and
$f_{m,n\pm1}$ to $$\matrix{
	T^{\lambda}_{m,n}f_{m,n} &=& T^{\lambda}_{m,n}f_{m,n} \hfill\cr
	T^{\lambda}_{m-1,n}f_{m-1,n} &=& T^{\lambda}_{m,n}f_{m-1,n}
		(1-\lambda a_{m-1,n})^{-1} \cr
	T^{\lambda}_{m+1,n}f_{m+1,n} &=& T^{\lambda}_{m,n}f_{m+1,n} \hfill\cr
	T^{\lambda}_{m,n-1}f_{m,n-1} &=& T^{\lambda}_{m,n}f_{m,n-1}
		(1-\lambda b_{m,n-1})^{-1} \cr
	T^{\lambda}_{m,n+1}f_{m,n+1} &=& T^{\lambda}_{m,n}f_{m,n+1} \hfill\cr
	}\Eqno\Label\MDeformation$$
i.e.\ any vertex star and its image under the $T^{\lambda}$ are M\"obius
equivalent.
In this sense, $\lambda\to T^{\lambda}_{m,n}f_{m,n}$ can be considered a
``M\"obius-deformation'' of $f$:

\Definition{~(M\"obius-deformation)}{Let $f:\Z^2\to\HP^1$ be a discrete
net; a family $(T^{\lambda})_{\lambda}$, $T^{\lambda}:\Z^2\to Gl(2,\H)$,
will be called a M\"obius-deformation for $f$ if any corresponding vertex
stars of $f=T^0f$ and $T^{\lambda}f$ are M\"obius equivalent.}

The integrability conditions for \TSystem\ directly lead to a discrete analog
for the ``general $C$-transform'' for an isothermic net $f:\Z^2\to\H$ \BPP:
\PLabel\GenCTransform

{\Lemma{~($C$-transformation)}{The system \TSystem\ is integrable for all
$\lambda\in\R$ if and only if $U=\partial_1{\fk F}^{\ast}$ and $V=\partial_2
{\fk F}^{\ast}$ with some ${\fk F}^{\ast}:\Z^2\to{\fk gl}(2,\H)$:
$$\matrix{
	{\fk F}^{\ast}_{m+1,n} &=& {\fk F}^{\ast}_{m,n} + U_{m,n} &,&
	{\fk F}^{\ast}_{m,n+1} &=& {\fk F}^{\ast}_{m,n} + V_{m,n} &.&\cr}
	\Eqno\Label\CSystem$$}
\Proof
The discrete version of the Maurer-Cartan equations for the system \TSystem\
reads $(1+\lambda U_{m,n})(1+\lambda V_{m+1,n})=(1+\lambda V_{m,n})(1+\lambda
U_{m,n+1})$; as $\varphi\circ f\equiv0$ this simplifies to $\lambda(U_{m,n}+
V_{m+1,n})=\lambda(V_{m,n}+U_{m,n+1})$, which is the integrability condition
for \CSystem\ when $\lambda=1$.
\Endproof}

\SubSection{The $C$-transformation}
Certain projections then give (as we will see a bit later) the well known
\BoPiB\ discrete analog of the Christoffel transformation in Euclidean space
(cf.\BPP):

{\Lemma{}{Let $v_0,v_{\infty}\in\H^2$ be homogeneous coordinates of
$0,\infty\in\HP^1$. Then, \CSystem\ is integrable if and only if there is a
discrete map ${\fk f}^{\ast}=\nu_{\infty}{\fk F}^{\ast}v_{\infty}:\Z^2\to\H$
such that $u=\partial_1{\fk f}^{\ast}$ and $v=\partial_2{\fk f}^{\ast}$,
and $$\matrix{
	(\partial_1{\fk f}^{\ast})_{m,n}(\partial_2{\fk f})_{m+1,n} &=&
	(\partial_2{\fk f}^{\ast})_{m,n}(\partial_1{\fk f})_{m,n+1}, \cr
	(\partial_1{\fk f})_{m,n}(\partial_2{\fk f}^{\ast})_{m+1,n} &=&
	(\partial_2{\fk f})_{m,n}(\partial_1{\fk f}^{\ast})_{m,n+1}. \cr
	}\Eqno\Label\DualRelations$$
Moreover, with
${\fk d}_{m,n}:=(\partial_1{\fk f})_{m,n}-(\partial_2{\fk f})_{m,n}
=(\partial_1{\fk f})_{m,n+1}-(\partial_2{\fk f})_{m+1,n}$, $$\matrix{
	{\fk d}_{m,n}
	(\partial_1{\fk f}^{\ast})_{m,n+1}(\partial_1{\fk f})_{m,n+1} &=&
	(\partial_1{\fk f})_{m,n}(\partial_1{\fk f}^{\ast})_{m,n}
		{\fk d}_{m,n},\cr
	{\fk d}_{m,n}
	(\partial_2{\fk f}^{\ast})_{m+1,n}(\partial_2{\fk f})_{m+1,n} &=&
	(\partial_2{\fk f})_{m,n}(\partial_2{\fk f}^{\ast})_{m,n}
		{\fk d}_{m,n};\cr
	{\fk d}_{m,n}
	(\partial_2{\fk f}^{\ast})_{m+1,n}(\partial_1{\fk f})_{m,n+1} &=&
	(\partial_2{\fk f})_{m,n}(\partial_1{\fk f}^{\ast})_{m,n}
		{\fk d}_{m,n},\cr
	{\fk d}_{m,n}
	(\partial_1{\fk f}^{\ast})_{m,n+1}(\partial_2{\fk f})_{m+1,n} &=&
	(\partial_1{\fk f})_{m,n}(\partial_2{\fk f}^{\ast})_{m,n}
		{\fk d}_{m,n}.\cr }\Eqno\Label\DualConsequences$$
	}\PLabel\CLemma
\Proof
A straightforward computation shows that the integrability condition for
\OldTSystem\ is equivalent to $u_{m,n}+v_{m+1,n}=v_{m,n}+u_{m,n+1}$,
i.e.\ $u=\partial_1{\fk f}^{\ast}$ and $v=\partial_2{\fk f}^{\ast}$
for some discrete map $f^{\ast}:\Z^2\to\H$, together with \DualRelations.
Since $\nu_{\infty}f=1=\varphi v_{\infty}$, $\nu_{\infty}Uv_{\infty}=u$
and $\nu_{\infty}Vv_{\infty}=v$ and hence ${\fk f}^{\ast}=\nu_{\infty}
{\fk F}^{\ast}v_{\infty}$.
%
%
Moreover, using \DualRelations, each of the equations \DualConsequences\
becomes the integrability condition for ${\fk f}^{\ast}$.
\Endproof}

Note that \DualRelations\ yields a discrete version of the equations
$d{\fk f}\wedge d{\fk f}^{\ast}=0=d{\fk f}^{\ast}\wedge d{\fk f}$
characterizing the Christoffel transform ${\fk f}^{\ast}:M^2\to\H$
in the smooth case \JePe.

{
Assuming that ${\fk f}^{\ast}$ is a discrete {\it net}, another useful
formula for the cross ratios of elementary quadrilaterals of $f$ can be
obtained as a direct consequence of \DualRelations:

\Lemma{}{$[f_{m,n},f_{m+1,n},f_{m+1,n+1},f_{m,n+1}]={\rm Re}\,q_{m,n}+
i|{\rm Im}\,q_{m,n}|$ where
$$q_{m,n}=[(\partial_1{\fk f}^{\ast})(\partial_1{\fk f})]_{m,n}
	[(\partial_2{\fk f}^{\ast})(\partial_2{\fk f})]_{m,n}^{-1}
	\Eqno\Label\CrossRelation$$
if ${\fk f}^{\ast}$ is a discrete net, satisfying \DualRelations\
with some stereographic projection ${\fk f}$ of $f$.}}

If the functions $a:=(\partial_1{\fk f}^{\ast})(\partial_1{\fk f}),
b:=(\partial_2{\fk f}^{\ast})(\partial_2{\fk f}):\Z^2\to\R$ are real
valued, then the first two equations of \DualConsequences\ show (the
real part of a quaternionic product is symmetric) that $a=a_m$ and
$b=b_n$ depend on $m$ resp.\ $n$ only.
Thus, \CrossRelation\ implies that $f:\Z^2\to\HP^1$ is a discrete isothermic
net; and, ${\fk f}^{\ast}$ is a Christoffel transform \BoPiB\ of the
stereographic projection ${\fk f}:\Z^2\to\H$ of $f$ --- that, for
symmetry reasons, is an isothermic net with $q^{\ast}=q$, too:

\Definition{~($C$-transformation)}{Let $f:\Z^2\to\HP^1$ be a discrete
isothermic net, $q_{m,n}={a_m\over b_n}$, ${\fk f}=(\nu_0 f)
(\nu_{\infty}f)^{-1}:\Z^2\to\H$ a stereographic projection of $f$.
Then, the discrete isothermic net ${\fk f}^{\ast}:\Z^2\to\H$ given
by $(\partial_1{\fk f}^{\ast})=a\,(\partial_1{\fk f})^{-1}$ and
$(\partial_2{\fk f}^{\ast})=b\,(\partial_2{\fk f})^{-1}$ is called
a $C$-transform (Christoffel transform) of ${\fk f}$.}

Note, that the $C$-transform ${\fk f}^{\ast}$ of ${\fk f}$ is unique
up to scale and translation in $\H$.

\SubSection{The $G$-transformation}
Obviously, the $C$-transformation depends on the Euclidean structure
of the ambient space; on the other hand, the notion of an isothermic
net is a conformal notion. Using this interplay of geometries, a new
transformation can be introduced --- as in the smooth case where this
transformation generalizes the classical Goursat transformation for
minimal surfaces \JMN:

\Definition{~($G$-transformation)}{Two isotermic nets ${\fk f},\tilde{\fk f}:
Z^2\to\H$ are said to be $G$-transforms (Goursat transforms) of each other
if their $C$-transforms are stereographic projections of the same isothermic
net $f:\Z^2\to\HP^1$.}

As the $C$-transformation is involutive, ${\fk f}^{\ast\ast}\simeq{\fk f}$
up to scaling and translation, the $G$-transformations of an isothermic
net ${\fk f}:\Z^2\to\H$ form a group.  Generically, there is a 4-parameter
family of $G$-transforms of a given isothermic net that are not congruent
to the original net: if ${\fk f}^{\ast}$ and $\tilde{\fk f}^{\ast}$ denote
$C$-transforms of different stereographic projections ${\fk f}=(\nu_0f)$ and
$\tilde{\fk f}=(\tilde\nu_0f)(\tilde\nu_{\infty}f)^{-1}$ of an isothermic net
$f:\Z^2\to\HP^1$ (as before, without loss of generality $\nu_{\infty}f\equiv1$),
then \Edges\ yields $$\matrix{
	&&(\partial_1\tilde{\fk f})_{m,n} \hfill&=&
		(\varphi_{m+1,n}\tilde{v}_{\infty})^{-1}
		(\partial_1{\fk f})_{m,n}
		(\tilde\nu_{\infty}f_{m,n})^{-1} \cr
	&\Leftrightarrow&
	(\partial_1\tilde{\fk f}^{\ast})_{m,n} &=&\hfill
		(\tilde\nu_{\infty}f_{m,n})
		(\partial_1{\fk f}^{\ast})_{m,n}
		(\varphi_{m+1,n}\tilde{v}_{\infty}), \cr
	&&(\partial_2\tilde{\fk f})_{m,n} \hfill&=&
		(\varphi_{m,n+1}\tilde{v}_{\infty})^{-1}
		(\partial_2{\fk f})_{m,n}
		(\tilde\nu_{\infty}f_{m,n})^{-1} \cr
	&\Leftrightarrow&
	(\partial_2\tilde{\fk f}^{\ast})_{m,n} &=&\hfill
		(\tilde\nu_{\infty}f_{m,n})
		(\partial_2{\fk f}^{\ast})_{m,n}
		(\varphi_{m,n+1}\tilde{v}_{\infty}). \cr
	}\Eqno\Label\GSystem$$
Thus, ${\fk f}^{\ast}$ and $\tilde{\fk f}^{\ast}$ will generically be
non-congruent if $\tilde{v}_{\infty}\H\neq v_{\infty}\H$ are different
points in $\HP^1$.

Clearly, any ``Christoffel pair'' ${\fk f},{\fk f}^{\ast}:\Z^2\to\H$ in
$\H$, consisting of a discrete isothermic net ${\fk f}$ and a Christoffel
transform ${\fk f}^{\ast}$ of ${\fk f}$, satisfies \DualRelations.
Thus, \CSystem\ is integrable for $f=v_0+v_{\infty}{\fk f}$ and
$\varphi=\nu_0-{\fk f}\nu_{\infty}$, $u=\partial_1{\fk f}^{\ast}$ and
$v=\partial_2{\fk f}^{\ast}$:
the ``general $C$-transform'' ${\fk F}^{\ast}:\Z^2\to{\fk gl}(2,\H)$ of an
isothermic net $f:\Z^2\to\HP^1$ encodes the $C$-transforms of all possible
stereographic projections of $f$ --- therefore, any two ``projections''
$\nu_{\infty}{\fk F}^{\ast}v_{\infty}$ and $\tilde\nu_{\infty}{\fk F}^{\ast}
\tilde{v}_{\infty}$ into Euclidean space are $G$-transforms of each other.
Choosing another stereographic projection $\tilde{\fk f}=(\tilde\nu_0f)
(\tilde\nu_{\infty}f)^{-1}$ of $f$ from the beginning, the $C$-transform
${\fk f}^{\ast}$ changes by a $G$-transformation --- \GSystem\ then shows
that $U$ and $V$ in \RForm, $$\matrix{
	\tilde U_{m,n} &=&
		[f(\tilde\nu_{\infty}f)^{-1}]_{m,n}
		(\partial_1\tilde{\fk f}^{\ast})_{m,n}
		[(\varphi\tilde{v}_{\infty})^{-1}\varphi]_{m+1,n}
		&=& U_{m,n} \cr
	\tilde V_{m,n} &=&
		[f(\tilde\nu_{\infty}f)^{-1}]_{m,n}
		(\partial_2\tilde{\fk f}^{\ast})_{m,n}
		[(\varphi\tilde{v}_{\infty})^{-1}\varphi]_{m,n+1}
		&=& V_{m,n}, \cr}$$
do not depend on the choice of stereographic projection ${\fk f}$.
Hence, ${\fk F}^{\ast}$ itself does not depend on the choice of stereographic
projection --- and, ``contains'' the $G$-equivalence class of all possible
$C$-transforms of an isothermic net $f:\Z^2\to\HP^1$ as ``off-diagonal
elements'' ${\fk f}^{\ast}=\nu_{\infty}{\fk F}^{\ast}v_{\infty}$:
\Definition{~($C$-transformation)}{Given an isothermic net $f:\Z^2\to\HP^1$,
a solution ${\fk F}^{\ast}:\Z^2\to{\fk gl}(2,\H)$ of \CSystem\ is called the
{\rm general} $C$-transform of $f$.}

\SubSection{The $T$-transformation}
With the same argument, a solution $T^{\lambda}$ of \TSystem\ does not depend
on the choice of stereographic projection.
Consequently, the following definition makes sense:

\Definition{~($T$-transformation)}{Let $f:\Z^2\to\HP^1$ be an isothermic net,
${\fk f}^{\ast}$ a $C$-transform of any stereographic projection ${\fk f}
=(\nu_0f)(\nu_{\infty}f)^{-1}$ of $f$; moreover, let $T^{\lambda}:\Z^2\to
Gl(2,\H)$ be a solution of \TSystem\ with $u=\partial_1{\fk f}^{\ast}$ and
$v=\partial_2{\fk f}^{\ast}$: $$\matrix{
	T^{\lambda}_{m+1,n} &=& T^{\lambda}_{m,n}(1+\lambda
		[f(\nu_{\infty}f)^{-1}]_{m,n}
		(\partial_1{\fk f}^{\ast})_{m,n}
		[(\varphi v_{\infty})^{-1}\varphi]_{m+1,n}), \cr
	T^{\lambda}_{m,n+1} &=& T^{\lambda}_{m,n}(1+\lambda
		[f(\nu_{\infty}f)^{-1}]_{m,n}
		(\partial_2{\fk f}^{\ast})_{m,n}
		[(\varphi v_{\infty})^{-1}\varphi]_{m,n+1}). \cr
	}\eqno{\rm(T)}$$
Then, $f^{\lambda}:=T^{\lambda}f:\Z^2\to\HP^1$ is called a $T$-transform
of $f$.}

As the $C$-transform ${\fk f}^{\ast}$ is unique up to scale and translation,
the family $(f^{\lambda})_{\lambda\in\R}$ of $T$-transforms of $f$ is unique
up to M\"obius transformations; fixing a factorization of the cross ratios,
$q_{m,n}={a_m\over b_n}$, each $f^{\lambda}$ is unique up to a M\"obius
transformation.

Generically, a $T$-transform $f^{\lambda}$ of an isothermic net $f$ is
a discrete net:

{\Lemma{}{Let $f:\Z^2\to\HP^1$ be an isothermic net, $q_{m,n}={a_m\over b_n}$.
Choosing an invertible initial condition, $T^{\lambda}_{0,0}\in Gl(2,\H)$, a
solution $T^{\lambda}$ of (T) stays invertible if and only if $\lambda a_m\neq
1\neq\lambda b_n$.}
\Proof
By \TFixpoints, $(1+\lambda U),(1+\lambda V):\Z^2\to Gl(2,\H)$ iff
$\lambda a_m\neq1\neq\lambda b_n$ since the regularity condition on
$f$ ($f$ is a discrete {\it net}) implies that $f_{m,n},f_{m+1,n}$
and $f_{m,n},f_{m,n+1}$ are (quaternionic) linearly independent.
\Endproof}

Moreover, if $f^{\lambda}$ is a discrete net, then it is isothermic:

{\Lemma{}{Let $f^{\lambda}$ be a non-degenerate $T$-transform of an isothermic
net $f$ with $q_{m,n}={a_m\over b_n}$; then $f^{\lambda}$ is isothermic with
$q^{\lambda}_{m,n}=q_{m,n}{1-\lambda b_n\over 1-\lambda a_m}={a_m\over
1-\lambda a_m}:{b_n\over 1-\lambda b_n}$.}
\Proof
$\varphi^{\lambda}=\varphi(T^{\lambda})^{-1}$ as $f^{\lambda}=T^{\lambda}f$.
Consequently, using \TFixpoints, $$\matrix{
	\hfill \varphi^{\lambda}_{m,n}f^{\lambda}_{m+1,n}
		&=& \hfill \varphi_{m,n}(1+\lambda U_{m,n})f_{m+1,n}
		&=& \varphi_{m,n}f_{m+1,n}, \cr
	\hfill \varphi^{\lambda}_{m+1,n+1}f^{\lambda}_{m+1,n}
		&=& \hfill \varphi_{m+1,n+1}(1+\lambda V_{m+1,n})^{-1}f_{m+1,n}
		&=& {\varphi_{m+1,n+1}f_{m+1,n}\over1-\lambda b_n}, \cr
	\hfill \varphi^{\lambda}_{m+1,n+1}f^{\lambda}_{m,n+1}
		&=& \hfill \varphi_{m+1,n+1}(1+\lambda U_{m,n+1})^{-1}f_{m,n+1}
		&=& {\varphi_{m+1,n+1}f_{m,n+1}\over1-\lambda a_m}, \cr
	\hfill \varphi^{\lambda}_{m,n}f^{\lambda}_{m,n+1}
		&=& \hfill \varphi_{m,n}(1+\lambda V_{m,n})f_{m,n+1}
		&=& \varphi_{m,n}f_{m,n+1}. \cr
	}$$ Thus, \CrossRatio\ yields the result
$q^{\lambda}_{m,n}=q_{m,n}{1-\lambda b_n\over1-\lambda a_m}$.
\Endproof}

This shows that the ``regular'' $T$-transformations act on the space of
(M\"obius equivalence classes of) isothermic nets. 
The next lemma shows that (for ``regular values'' of $\lambda$) the
$T$-transformations act like a ``1-parameter group'' on the (1-dimensional)
orbit $(f^{\lambda})_{\lambda}$ of a fixed isothermic net:

{\Theorem{}{Let $f=f^0:\Z^2\to\HP^1$ be an isothermic net and $(T^{\lambda})_%
{\lambda\in\R}$ its family of $T$-transformations, $T^{\lambda}_{0,0}=id$ for
all $\lambda\in\R$.
If $T^{\lambda_1},T^{\lambda_2}:\Z^2\to Gl(2,\H)$ are non-degenerate, then
$T^{\lambda_1+\lambda_2}=T^{\lambda_2}T^{\lambda_1}$.
In particular, $(T^{\lambda})^{-1}=T^{-\lambda}$.}
\Proof
Using \TFixpoints,
$T^{\lambda_1}_{m+1,n}=T^{\lambda_1}_{m,n}(1+\lambda_1U_{m,n})$ and
$T^{\lambda_1}_{m,n+1}=T^{\lambda_1}_{m,n}(1+\lambda_1V_{m,n})$,
$$\matrix{
	(1+\lambda_2U^{\lambda_1}_{m,n})T^{\lambda_1}_{m,n}
		(1+\lambda_1U_{m,n})f_{m,n} \hfill&=&
		f^{\lambda_1}_{m,n}(1-(\lambda_1+\lambda_2)a_m), \hfill\cr
	(1+\lambda_2U^{\lambda_1}_{m,n})T^{\lambda_1}_{m,n}
		(1+\lambda_1U_{m,n})f_{m+1,n} \hfill&=&
		f^{\lambda_1}_{m+1,n}, \hfill\cr
	(1+\lambda_2V^{\lambda_1}_{m,n})T^{\lambda_1}_{m,n}
		(1+\lambda_1V_{m,n})f_{m,n} \hfill&=&
		f^{\lambda_1}_{m,n}(1-(\lambda_1+\lambda_2)b_n), \hfill\cr
	(1+\lambda_2V^{\lambda_1}_{m,n})T^{\lambda_1}_{m,n}
		(1+\lambda_1V_{m,n})f_{m,n+1} \hfill&=&
		f^{\lambda_1}_{m,n+1}. \hfill\cr
	}$$
Hence,
$$\matrix{
	(T^{\lambda_1})^{-1}(1+\lambda_2U^{\lambda_1})T^{\lambda_1}
		(1+\lambda_1U) &=& 1+(\lambda_1+\lambda_2)U, \cr
	(T^{\lambda_1})^{-1}(1+\lambda_2V^{\lambda_1})T^{\lambda_1}
		(1+\lambda_1V) &=& 1+(\lambda_1+\lambda_2)V. \cr
	}$$
Therefore, $T^{\lambda_2}T^{\lambda_1}$ and $T^{\lambda_1+\lambda_2}$
satisfy the same difference equations; since they coincide at $(m,n)=(0,0)$,
$T^{\lambda_2}T^{\lambda_1}=T^{\lambda_1+\lambda_2}$.
In particular, $(T^{\lambda})^{-1}=T^{-\lambda}$.
\Endproof}

Finally, we obtain some type of discrete version of the non-rigidity result
for isothermic surfaces in 3-dimensional M\"obius geometry:

{\Lemma{}{Let $f:\Z^2\to S^3$ be a discrete net in a 3-sphere
$S^3\subset\HP^1$; then, a solution $T^{\lambda}:\Z^2\to M\ddot{o}b(S^3)$
of \TSystem\ takes --- up to a fixed M\"obius transformation --- values
in the M\"obius group of $S^3$ if and only if $f$ is isothermic.}
\Proof
Choose homogeneous coordinates $s\in\R^6_1$ for $S^3$, and
let $a:=(\partial_1{\fk f}^{\ast})(\partial_1{\fk f})$ and
$b:=(\partial_2{\fk f}^{\ast})(\partial_2{\fk f})$; then, \TFixpoints\
yields $$\matrix{
	((1+\lambda U_{m,n})^{-1}\cdot s)(f_{m,n},f_{m,n}) \hfill&=&
	((1+\lambda U_{m,n})^{-1}\cdot s)(f_{m+1,n},f_{m+1,n}) \hfill\cr
		&=& 0 \hfill {\rm\ and} \hfill\cr
	((1+\lambda U_{m,n})^{-1}\cdot s)(f_{m+1,n},f_{m,n}) &=&
		s(f_{m+1,n},f_{m,n})(1-\lambda a_{m,n})^{-1}; \hfill\cr
	((1+\lambda V_{m,n})^{-1}\cdot s)(f_{m,n},f_{m,n}) \hfill&=&
	((1+\lambda V_{m,n})^{-1}\cdot s)(f_{m,n+1},f_{m,n+1}) \hfill\cr
		&=& 0 \hfill {\rm\ and} \hfill\cr
	((1+\lambda V_{m,n})^{-1}\cdot s)(f_{m,n+1},f_{m,n}) &=&
		s(f_{m,n+1},f_{m,n})(1-\lambda b_{m,n})^{-1}. \hfill\cr
	}$$
Thus, $(1+\lambda U),(1+\lambda V):Z^2\to M\ddot{o}b(S^3)$ preserve the sphere
$S^3\cong\R s$ if and only if $a,b:\Z^2\to\R$ are real valued.
But, we already deduced it as a consequence of \DualConsequences, that this
is equivalent to $f$ being isothermic with $q_{m,n}={a_m\over b_n}$ where
the functions $a,b:\Z\to\R$ depend on one variable only.
\Endproof}

\CatenoidCousins

The following consequence will prove crucial for our approach to a
discrete analog for surfaces of constant mean curvature 1 in hyperbolic
space:

\Corollary{}{If $f:\Z^2\to S^3$ is an isothermic net, then (up to a
M\"obius transformation) its $T$-transforms $f^{\lambda}:\Z^2\to S^3$
take values in $S^3$, too.\newline
If $f:\Z^2\to S^2$, then $f^{\lambda}:\Z^2\to S^2$.}

\SubSection{The $D$-transformation}
Finally, we formulate the discrete analog of the Darboux transformation
in the conformal setting:

\Definition{~($D$-transformation)}{A discrete net $\hat{f}:\Z^2\to\HP^1$
is called a $D$-transform (Darboux transform) of an isothermic net
$f:\Z^2\to\HP^1$, $\hat{f}=D_{\lambda}f$, if $T^{\lambda}\hat{f}
=const\in\HP^1$ is a fixed point.}

Obviously, for fixed $\lambda$ and fixed ``initial point'' $\hat{f}_{0,0}
\in\HP^1$ there is a unique Darboux transform
$(T^{\lambda})^{-1}(T^{\lambda}_{0,0}\hat{f}_{0,0}):\Z^2\to\HP^1$ of an
isothermic net $f$.  If $f$ is an isothermic net on a 3-sphere $S^3$, then
$T^{\lambda}:Z^2\to M\ddot{o}b(S^3)$ so that any $D$-transform with initial
point $\hat{f}_{0,0}\in S^3$ stays on $S^3$.
Thus, a fixed isothermic net $f$ allows a $\infty^{1+4}$ $D$-transforms in
$\HP^1$, or $\infty^{1+3}$ $D$-transforms in $S^3$.

Actually, the above definition for the $D$-transformation agrees with the
one given in \HHP, where we used a discrete version of the Riccati equation
characterizing smooth Darboux transforms in Euclidean 4-space:

{\Lemma{}{A discrete net $\hat{f}:\Z^2\to\HP^1$ is a Darboux transform
of an isothermic net $f:\Z^2\to\HP^1$ if and only if their stereographic
projections ${\fk f},\hat{\fk f}:Z^2\to\H$ satisfy the Riccati type system
$$\matrix{
	(\partial_1\hat{\fk f})_{m,n} &=& \lambda
		(\hat{\fk f}-{\fk f})_{m,n}
		(\partial_1{\fk f}^{\ast})_{m,n}
		(\hat{\fk f}-{\fk f})_{m+1,n}, \cr
	(\partial_2\hat{\fk f})_{m,n} &=& \lambda
		(\hat{\fk f}-{\fk f})_{m,n}
		(\partial_2{\fk f}^{\ast})_{m,n}
		(\hat{\fk f}-{\fk f})_{m,n+1}. \cr
	}\Eqno\Label\Riccati$$}
\Proof
$\hat{f}:\Z^2\to\HP^1$ being a $D$-transform of an isothermic net
$f:\Z^2\to\HP^1$ is equivalent to $(1+\lambda U_{m,n})\hat{f}_{m+1,n}
\parallel\hat{f}_{m,n}\parallel(1+\lambda V_{m,n})\hat{f}_{m,n+1}$
for all $m,n$. As an equivalent formulation, we obtain $$\matrix{
	0 &=& \hat\varphi_{m,n}(1+\lambda U_{m,n})\hat{f}_{m+1,n} \hfill\cr
	&=& (\partial_1\hat{\fk f})_{m,n} + \lambda
		({\fk f}-\hat{\fk f})_{m,n}
		(\partial_1{\fk f}^{\ast})_{m,n}
		(\hat{\fk f}-{\fk f})_{m+1,n} \cr
	0 &=& \hat\varphi_{m,n}(1+\lambda V_{m,n})\hat{f}_{m,n+1} \hfill\cr
	&=& (\partial_2\hat{\fk f})_{m,n} + \lambda
		({\fk f}-\hat{\fk f})_{m,n}
		(\partial_2{\fk f}^{\ast})_{m,n}
		(\hat{\fk f}-{\fk f})_{m,n+1} \cr
	}$$ from \Edges\ --- remember that we agreed on having
$\nu_{\infty}\hat{f}\equiv1\equiv\hat{\varphi}v_{\infty}$.
\Endproof}

This lemma allows to carry over the results given in \HHP\ to the more
general setting\Footnote{In \HHP, we restricted to the special case $a_m=+1$,
$b_n=-1$ of the discrete analog of {\it conformal} curvature line coordinates
of isothermic surfaces in $\H\cong\R^4$.} of isothermic nets $f:\Z^2\to\HP^1$
``in the wide sense''.
As in \HHP, \Riccati\ can be reformulated to obtain cross ratio conditions
for the quadrilaterals formed by corresponding edges on $f$ and $\hat{f}$:
$$\matrix{
	\lambda a_m &=& [f_{m,n},f_{m+1,n},\hat{f}_{m+1,n},\hat{f}_{m,n}], \cr
	\lambda b_n &=& [f_{m,n},f_{m,n+1},\hat{f}_{m,n+1},\hat{f}_{m,n}]. \cr
	}\Eqno\Label\CRiccati$$
With this observation, several key results on the Darboux transformation can
be reduced (cf.\HHP) to the

\Lemma{~(Hexahedron lemma)}{Let $p_1,\dots,p_4\in\HP^1$ be four concircular
points, i.e.\ $[p_1,p_2,p_3,p_4]=:\mu\in\R$, and $\lambda\in\R$. Then, to any
$q_1\in\HP^1$ there exist unique points $q_2,q_3,q_4\in\HP^1$ such that
$$\matrix{
	[q_1,q_2,q_3,q_4] &=& [p_1,p_2,p_3,p_4] &=& \mu \cr
	[q_1,q_2,p_2,p_1] &=& [q_3,q_4,p_4,p_3] &=& \mu\lambda \cr
	[q_2,q_3,p_3,p_2] &=& [q_4,q_1,p_1,p_4] &=& \lambda \cr
	}$$
Moreover, the eight points $p_i,q_i$ lie on a 2-sphere.}

As a direct consequence (cf.\HHP,\BoJe),

\Theorem{}{A $D$-transform $\hat{f}:Z^2\to\HP^1$ of an isothermic
net $f:\Z^2\to\HP^1$ is an isothermic net with the same cross ratios
$\hat{q}_{m,n}=q_{m,n}$;
moreover, $f$ and $\hat{f}$ envelope a discrete Ribaucour congruence
of 2-spheres $s_1\wedge s_2:\Z^2\to\Lambda^2(\R^6_1)$, i.e.\ any two
corresponding elementary quadrilaterals of $f$ and $\hat{f}$ lie on
a two sphere $S^2_{m,n}=(S_1)_{m,n}\cap(S_2)_{m,n}\subset\HP^1$.}

As a consequence, the Darboux transformation is involutive: since
\CRiccati\ is symmetric in $f$ and $\hat{f}$, $f$ is a Darboux
transform of $\hat{f}$ as soon as $\hat{f}$ is known to be isothermic
with the same cross ratio as $f$.
Thus, any isothermic net together with a $D$-transform $\hat{f}$
of $f$ form a ``Darboux pair'' of isothermic nets.

Permutability of different $D$-transforms is, again, a direct consequence
of the hexahedron lemma (cf.\HHP):

\Theorem{}{Given two $D$-transforms $\hat{f}_1=D_{\lambda_1}f$ and
$\hat{f}_1=D_{\lambda_1}f$ of an isothermic net $f:\Z^2\to\HP^1$, there
is an isothermic net $\hat{f}=D_{\lambda_2}\hat{f}_1=D_{\lambda_1}\hat{f}_2$:
in this sense, $D_{\lambda_1}D_{\lambda_2}=D_{\lambda_2}D_{\lambda_1}$;
$\hat{f}$ is given by the relation $[f,\hat{f}_2,\hat{f},\hat{f}_1]
\equiv{\lambda_1\over\lambda_2}$.}

\SubSection{Permutability}
theorems also control the interrelations of different transformations of
isothermic nets: all the various transformations are related by certain
``commutation relations''.
We start with the permutability of the $C$- and $D$-transformations in
Euclidean space, giving a proof similar to the one in \HHP:

{\Theorem{}{The $C$-transforms of a $D_{\lambda}$-pair ${\fk f},\hat{\fk f}:
\Z^2\to\H$ of isothermic nets form (if properly scaled and positioned) a
$D_{\lambda}$-pair of isothermic nets: $D_{\lambda}C=CD_{\lambda}$}
\Proof
Let ${\fk f}^{\ast}$ be a $C$-transform of ${\fk f}$, and $\tilde{\fk f}
:={\fk f}^{\ast}+{1\over\lambda}(\hat{\fk f}-{\fk f})^{-1}$.
Then, \Riccati\ yields $$\matrix{
	(\partial_1\tilde{\fk f})_{m,n} &=& \lambda
		[{1\over\lambda}(\hat{\fk f}-{\fk f})^{-1}]_{m,n}
		(\partial_1{\fk f})_{m,n}
		[{1\over\lambda}(\hat{\fk f}-{\fk f})^{-1}]_{m+1,n}, \cr
	(\partial_2\tilde{\fk f})_{m,n} &=& \lambda
		[{1\over\lambda}(\hat{\fk f}-{\fk f})^{-1}]_{m,n}
		(\partial_2{\fk f})_{m,n}
		[{1\over\lambda}(\hat{\fk f}-{\fk f})^{-1}]_{m,n+1}, \cr
	}$$
i.e.\ $\tilde{\fk f}$ is a $D_{\lambda}$-transform of ${\fk f}^{\ast}$.
On the other hand, using the symmetries of the cross ratio, $
	(\partial_1\tilde{\fk f})(\partial_1\hat{\fk f}) =
	(\partial_1{\fk f})(\partial_1{\fk f}^{\ast})$ and $
	(\partial_2\tilde{\fk f})(\partial_2\hat{\fk f}) =
	(\partial_2{\fk f})(\partial_2{\fk f}^{\ast})
	$ such that $\tilde{\fk f}=\hat{\fk f}^{\ast}$ also is a
$C$-transform of $\hat{\fk f}$.
\Endproof}

This permutability theorem can obviously be read the other way: given a
Christoffel pair ${\fk f},{\fk f}^{\ast}:\Z^2\to\H$ of isothermic nets,
then any suitably attuned Darboux transforms form a Christoffel pair, again.

The following gives the effect of the $T$-transformation on Christoffel pairs:

{\Lemma{}{Let ${\fk f},{\fk f}^{\ast}:\Z^2\to\H$ be a Christoffel-pair
of isothermic nets; $F^{\lambda}=T^{\lambda}F^0$ a solution of the
system \OldTSystem. Then, $F^{\lambda}v_{\infty}=(T^{\ast})^{\lambda}
(v_0+v_{\infty} {\fk f}^{\ast})$ is a $T$-transform of ${\fk f}^{\ast}
\simeq(v_0+v_{\infty}{\fk f}^{\ast})$ where $(T^{\ast})^{\lambda}=T^{\lambda}
(v_{\infty}\nu_{\infty}+\lambda(v_0+v_{\infty}{\fk f})(\nu_0-{\fk f}^{\ast}
\nu_{\infty}))$.}
\Proof
It is straightforward to check that $F^{\lambda}(v_{\infty}\nu_{\infty}
+\lambda v_0\nu_0)$ satisfies \OldTSystem\ with the roles of ${\fk f}$
and ${\fk f}^{\ast}$ interchanged; hence, $F^{\lambda}(v_{\infty}\nu_{\infty}
+\lambda v_0\nu_0)v_0=F^{\lambda}v_{\infty}$ is a $T$-transform of
$(v_0+v_{\infty}{\fk f}^{\ast})$.  With the Euclidean frame
$id+v_{\infty}{\fk f}^{\ast}\nu_{\infty}$ of ${\fk f}^{\ast}$,
$(T^{\ast})^{\lambda}=F^{\lambda}(v_{\infty}\nu_{\infty}+\lambda v_0\nu_0)
(id-v_{\infty}{\fk f}^{\ast}\nu_{\infty})=T^{\lambda}(v_{\infty}\nu_{\infty}
+\lambda(v_0+v_{\infty}{\fk f})(\nu_0-{\fk f}^{\ast}\nu_{\infty}))$
since, up to M\"obius transformation, the solutions of \OldTSystem\
and \TSystem\ are related by the Euclidean frame of the underlying
isothermic surface.
\Endproof}

{\Corollary{}{The $T$ transforms $T^{\lambda}(v_0+v_{\infty}{\fk f})$ and
$(T^{\ast})^{\lambda}(v_0+v_{\infty}{\fk f}^{\ast})$ of a Christoffel pair
${\fk f},{\fk f}^{\ast}:\Z^2\to\H$ form (if properly positioned) a Darboux
pair: $T^{\lambda}C=D_{-\lambda}T^{\lambda}$.}
\Proof
From the previous lemma, $(T^{\ast})^{\lambda}(v_0+v_{\infty}{\fk f}^{\ast})
=T^{\lambda}v_{\infty}$ with the ``proper positioning'' $(T^{\ast})^{\lambda}
=T^{\lambda}(v_{\infty}\nu_{\infty}+\lambda(v_0+v_{\infty}{\fk f})
(\nu_0-{\fk f}^{\ast}\nu_{\infty}))$ such that the claim follows since
$(T^{\lambda})^{-1}=T^{-\lambda}$.
\Endproof}

The converse of this statement is true, too: any Darboux pair of isothermic
nets comes from a Christoffel pair via $T$-transformation:

{\Theorem{}{The $T$-transforms $T^{\lambda}f$ and $\hat{T}^{\lambda}
\hat{f}$ of a $D_{\lambda}$-pair $f,\hat{f}:\Z^2\to\HP^1$ form (after
proper stereographic projection) a Christoffel pair:
$CT^{\lambda}=T^{\lambda}D_{\lambda}$.}
\Proof
Choose homogeneous coordinates $v_0,v_{\infty}\in\H^2$ of two points in
$\HP^1$.  As $\hat{f}$ is a $D_{\lambda}$-transform of $f$, without loss
of generality\Footnote{$T^{\lambda}\hat{f}\equiv const\in\HP^1$: thus, add
a M\"obius transformation to $T^{\lambda}$ to obtain $T^{\lambda}\hat{f}
\equiv v_{\infty}\in\HP^1$ and rescale $\hat{f}\to\hat{f}(\nu_0T^{\lambda}
\hat{f})^{-1}$; then, rescale $f\to f(\nu_{\infty}T^{\lambda}f)^{-1}$.}
$T^{\lambda}\hat{f}\equiv v_{\infty}\in\H^2$, and
$T^{\lambda}f=v_0+v_{\infty}{\fk f}$ with ${\fk f}:\Z^2\to\H$.
In this setting, $f=T^{-\lambda}(v_0+v_{\infty}{\fk f})$ and
$\hat{f}=T^{-\lambda}v_{\infty}$ with the $T$-transform $T^{-\lambda}:
\Z^2\to PGl(2,\H)$ associated to $v_0+v_{\infty} {\fk f}$ and its
$C$-transform ${\fk f}^{\ast}\simeq v_0+v_{\infty}{\fk f}^{\ast}$.
But,
$T^{-\lambda}v_{\infty}=(T^{\ast})^{-\lambda}(v_0+v_{\infty}{\fk f}^{\ast})$
is the $T^{-\lambda}$-transform of ${\fk f}^{\ast}$; thus, up to M\"obius
transformation $\hat{T}^{\lambda}=((T^{\ast})^{-\lambda})^{-1}$, and
hence $\hat{T}^{\lambda}\hat{f}=v_0+v_{\infty}{\fk f}^{\ast}$.
Consequently, with this positioning of $\hat{T}^{\lambda}\hat{f}$,
both nets, $T^{\lambda}f$ and $\hat{T}^{\lambda}\hat{f}$, simultanously
project to a $C$-pair.
\Endproof}

As an immediate consequence of the previous theorem, we see that the
``difference'' between two $D$-transforms can be ``measured'' by the
$G$-transformation:

\Corollary{}{The $T$-transforms $\hat{T}_1^{\lambda}\hat{f}_1$ and
$\hat{T}_2^{\lambda}\hat{f}_2$ of two (different) $D_{\lambda}$-transforms
$\hat{f}_1,\hat{f}_2:\Z^2\to\HP^1$ of an isothermic net $f:\Z^2\to\HP^1$
are (after proper stereographic projections) $G$-transforms of each other.}

Another consequence of the interrelations of $D$- and $C$-pairs under the
$T$-trans\-formation stated above is the ``generic'' behaviour of the Darboux
pairs under the $T$-transformation:

{\Corollary{}{The $T$-transforms $T^{\mu}f$ and $\hat{T}^{\mu}\hat{f}$,
$\mu\neq\lambda$, of a $D_{\lambda}$-pair of isothermic nets, $f,\hat{f}:
\Z^2\to\HP^1$, form (if properly positioned) a $D_{\lambda-\mu}$-pair:
in this sense, $T^{\mu}D_{\lambda}=D_{\lambda-\mu}T^{\mu}$.  Moreover,
$\hat{T}^{\mu}=T^{\mu}(1-{\mu\over\lambda}f(\hat\varphi f)^{-1}\hat\varphi)$.}
\Proof
By the previous theorem, $f=T^{-\lambda}(v_0+v_{\infty}{\fk f})$ and
$\hat{f}=(T^{\ast})^{-\lambda}(v_0+v_{\infty}{\fk f}^{\ast})$ for a
suitable $C$-pair ${\fk f},{\fk f}^{\ast}:\Z^2\to\H$; with
$T^{\mu}T^{-\lambda}=T^{\mu-\lambda}$, this implies the first statement.
The second statement follows by a straightforward computation:
$(T^{\ast})^{\mu-\lambda}((T^{\ast})^{\lambda})^{-1}
=T^{\mu}(1-{\mu\over\lambda}T^{-\lambda}(v_0+v_{\infty}{\fk f})\nu_{\infty}
(T^{-\lambda})^{-1})$.
\Endproof}

Combining the above permutability theorems on the effect of the
$D$- and $T$-transformations on Christoffel pairs, a larger
permutability scheme is obtained:
{\Theorem{}{The $T^{\lambda}$-transformation preserves the
$C$-$D_{\lambda}$-permutability: \smallskip\Figure{$$\matrix{
	\matrix{
	{\fk f} &\RArrow[20pt]_{C}^{}& {\fk f}^{\ast} \cr
	\DArrow[20pt]{}{D_{\lambda}} && \DArrow[20pt]{D_{\lambda}}{} \cr
	\hat{\fk f} &\RArrow[20pt]_{}^{C}& \hat{\fk f}^{\ast} \cr
	}&\RArrow[50pt]_{}^{T^{\lambda}}&\matrix{
	T^{\lambda}(v_0+v_{\infty}{\fk f})
		& \RArrow[40pt]_{D_{-\lambda}}^{} &
	(T^{\ast})^{\lambda}(v_0+v_{\infty}{\fk f}^{\ast}) \cr
	\DArrow[20pt]{}{C} && \DArrow[20pt]{C}{} \cr
	\hat{T}^{\lambda}(v_0+v_{\infty}\hat{\fk f})
		& \RArrow[40pt]_{}^{D_{-\lambda}} &
	(\hat{T}^{\ast})^{\lambda}(v_0+v_{\infty}\hat{\fk f}^{\ast}) \cr
	}\cr }$$}{A permutability theorem}\Label\PermThm}
\Proof
Given four isothermic nets ${\fk f},\hat{\fk f},\hat{\fk f}^{\ast},
{\fk f}^{\ast}:\Z^2\to\H$ that form two $D_{\lambda}$-pairs and two
$C$-pairs according to the $C$-$D_{\lambda}$-permutability theorem,
two $D_{-\lambda}$-pairs $T^{\lambda}(v_0+v_{\infty}{\fk f}),
(T^{\ast})^{\lambda}(v_0+v_{\infty}{\fk f}^{\ast})$ and
$\hat{T}^{\lambda}(v_0+v_{\infty}\hat{\fk f}),
(\hat{T}^{\ast})^{\lambda}(v_0+v_{\infty}\hat{\fk f}^{\ast})$
are obtained via suitably adjusted $T$-transformations:
$$\matrix{(T^{\ast})^{\lambda} &=& T^{\lambda}(v_{\infty}\nu_{\infty}
+\lambda(v_0+v_{\infty}{\fk f})(\nu_0-{\fk f}^{\ast}\nu_{\infty})),\cr
(\hat{T}^{\ast})^{\lambda} &=& \hat{T}^{\lambda}(v_{\infty}\nu_{\infty}
+\lambda(v_0+v_{\infty}\hat{\fk f})(\nu_0-\hat{\fk f}^{\ast}\nu_{\infty})).\cr
}$$
On the other hand, $T^{\lambda}(v_0+v_{\infty}\hat{\fk f})=(T^{\ast})^{\lambda}
(v_0+v_{\infty}\hat{\fk f}^{\ast})\in\HP^1$ define the same point
in $\HP^1$ since $(\hat{\fk f}^{\ast}-{\fk f}^{\ast})={1\over\lambda}
(\hat{\fk f}-{\fk f})^{-1}$; hence, $T^{\lambda}(v_0+v_{\infty}{\fk f}),
\hat{T}^{\lambda}(v_0+v_{\infty}\hat{\fk f})$ and
$(T^{\ast})^{\lambda}(v_0+v_{\infty}{\fk f}^{\ast}),
(\hat{T}^{\ast})^{\lambda}(v_0+v_{\infty}\hat{\fk f}^{\ast})$
simultanously project to $C$-pairs, when, additionally, $T^{\lambda}$
and $\hat{T}^{\lambda}$ (resp.\ $(T^{\ast})^{\lambda}$ and
$(\hat{T}^{\ast})^{\lambda}$) are ``properly adjusted'', as in the
proof of the $(CT^{\lambda}=T^{\lambda}D_{\lambda})$-permutability theorem.
\Endproof}


Having established this transformation theory for discrete isothermic nets,
we are prepared to study
\Section{Horospherical nets in hyperbolic space}
as a discrete  analog of surfaces of constant mean curvature 1 --- the
mean curvature of horospheres --- in hyperbolic space. Smooth cmc-1 surfaces
in hyperbolic space appear as $T$-transforms\Footnote{Recall that, in case of
constant mean curvature surfaces in space forms, the $T$-transformation becomes
transformation via Lawson correspondence; for minimal surfaces in Euclidean
space, the $T$-trans\-formation provides the Umehara-Yamada perturbation.}
(cousins) of minimal surfaces in Euclidean space, and, on the other hand,
as $D$-transforms of their hyperbolic Gauss maps \JMN.
These facts suggest two different ansatzes to define discrete horospherical
nets: as $T$-transforms of discrete minimal nets, or as $D$-transforms of
discrete isothermic nets in a 2-sphere.
Both ansatzes will turn out to give the same class of discrete nets.

First, recall the definition of discrete minimal nets \BoPiA:

\Definition{~(Minimal nets)}{An isothermic net ${\fk f}:Z^2\to{\rm Im}\H$
is called minimal if it is a $C$-transform ${\fk f}={\fk n}^{\ast}$ of an
isothermic net ${\fk n}:Z^2\to S^2\subset{\rm Im}\H$ into the 2-sphere;
${\fk n}$ is then called the Gauss map of the discrete minimal net.}

As an immediate consequence of this definition, discrete minimal nets can be
obtained from complex valued isothermic nets (``discrete holomorphic nets''):
the stereographic projection $i(i+{\fk g}j)(i-{\fk g}j)^{-1}:\Z^2\to S^2$ of
any discrete isothermic net ${\fk g}j:\Z^2\to\C j\cong\C$ appears as a Gauss
map (or, $C$-transform) of a discrete minimal net --- that, therefore, can
be (re-) constructed from ${\fk g}:\Z^2\to\C$.
A discrete analog of the classical Weierstrass representation for minimal
surfaces is obtained via the Goursat transformation: given an isothermic net
${\fk h}:\Z^2\to\C$, a minimal net ${\fk f}:\Z^2\to{\rm Im}\H$ is obtained as
the $C$-transform of the stereographic projection of\Footnote{Note, that we are
working with ``principal coordinates'' so that the minimal surface is already
determined by just {\it one} holomorphic function.} ${\fk g}={\fk h}^{\ast}$
--- with $\tilde v_0=(v_0-v_{\infty}i){i\over\sqrt{2}}$ and $\tilde v_{\infty}
=(v_0+v_{\infty}i){1\over\sqrt{2}}$, \GSystem\ yields a discrete version of
the (quaternionic: cf.\JMN) Weierstrass formula \BoPiB: $$\matrix{
	(\partial_1{\fk f})_{m,n} &=& {1\over2}(i-{\fk g}j)_{m,n}
		j(\partial_1{\fk h})_{m,n}(i-{\fk g}j)_{m+1,n}, \cr
	(\partial_2{\fk f})_{m,n} &=& {1\over2}(i-{\fk g}j)_{m,n}
		j(\partial_2{\fk h})_{m,n}(i-{\fk g}j)_{m,n+1}. \cr }$$

Another direct consequence of this definition for discrete minimal nets is
the promised equivalence of our two ansatzes to define a discrete analog of
constant mean curvature 1 surfaces in hyperbolic space:

{\Lemma{}{An isothermic net is a $T$-transform of a discrete minimal
net iff it is a $D$-transform of an isothermic net on the 2-sphere.}
\Proof
This is a consequence of two permutability theorems:
a discrete minimal net ${\fk f}:\Z^2\to{\rm Im}\H$ and its Gauss map
${\fk n}:\Z^2\to S^2$ form a $C$-pair, thus their $T^{-\lambda}$-transforms
$T^{-\lambda}(v_0+v_{\infty}{\fk f}):\Z^2\to S^3$ and $(T^{\ast})^{-\lambda}
(v_0+v_{\infty}{\fk n}):\Z^2\to S^2$ form (if properly positioned) a
$D_{\lambda}$-pair, $T^{-\lambda}C=D_{\lambda}T^{-\lambda}$;
on the other hand, the $T^{\lambda}$ transforms of a $D_{\lambda}$-pair
$f:\Z^2\to S^3$ and $n:\Z^2\to S^2$ form (after proper stereographic
projection) a $C$-pair, $T^{\lambda}D_{\lambda}=CT^{\lambda}$, and
$T^{\lambda}n:\Z^2\to S^2$ takes values in some 2-sphere while
$T^{\lambda}f:\Z^2\to{\rm Im}\H$.
\Endproof}

Thus, similar to the definition of minimal nets, we define a discrete analog
for constant mean curvature 1 surfaces in hyperbolic space via coupling with
their hyperbolic Gauss maps. This determines $H^3\subset S^3$ as one of the
connected components of $S^3\setminus S^2$ of the complement of its infinity
boundary $S^2\cong\partial H^3$:

\Definition{~(Horospherical nets)}{An isothermic net $f:\Z^2\to S^3\setminus
S^2$ is called horospherical if it is a $D$-transform of an isothermic net
$n:\Z^2\to S^2\cong\partial H^3$;
this net, $n$, is then called the hyperbolic Gauss map of $f$.}

This definition directly leads to a representation in terms of ``discrete
holomorphic'' data: given a discrete isothermic net ${\fk g}:\Z^2\to\C$
together with its $C$-transform ${\fk h}:={\fk g}^{\ast}:\Z^2\to\C$ (i.e.\
$({\fk g}j)^{\ast}=-j{\fk h}$), any $D$-transform of $n=v_0+v_{\infty}{\fk g}j$
comes from a solution $T^{\lambda}:\Z^2\to Gl(2,\H)$ of (T),
$$\matrix{
	T^{\lambda}_{m+1,n}
	&=& T^{\lambda}_{m,n}(1-\lambda[v_0+v_{\infty}{\fk g}_{m,n}j]
		j(\partial_1{\fk h})_{m,n}
		[\nu_0-{\fk g}_{m+1,n}j\nu_{\infty}]) \cr
	&=& T^{\lambda}_{m,n}J^{-1}(1+\lambda[v_0+v_{\infty}{\fk g}_{m,n}]
		(\partial_1{\fk h})_{m,n}
		[\nu_0-{\fk g}_{m+1,n}\nu_{\infty}]) J, \cr
	T^{\lambda}_{m,n+1}
	&=& T^{\lambda}_{m,n}(1+\lambda[v_0+v_{\infty}{\fk g}_{m,n}j]
		j(\partial_2{\fk h})_{m,n}
		[\nu_0-{\fk g}_{m,n+1}j\nu_{\infty}]) \cr
	&=& T^{\lambda}_{m,n}J^{-1}(1-\lambda[v_0+v_{\infty}{\fk g}_{m,n}]
		(\partial_2{\fk h})_{m,n}
		[\nu_0-{\fk g}_{m,n+1}\nu_{\infty}]) J, \cr
	}\Eqno\Label\WSystem$$
via $D_{\lambda}n=(T^{-\lambda})^{-1}p_0$ where $p_0\in\HP^1$ is a (constant)
point.
Here, we denoted $J:=(v_{\infty}\nu_0+v_0j\nu_{\infty})$: obviously, \WSystem\
is a {\it complex} system for $T^{\lambda}J^{-1}$ --- indeed, it is a discrete
version of the system arising in Bryant's Weierstrass type representation for
constant mean curvature 1 surfaces in hyperbolic space \UmYa. Thus, $f^{\sharp}
=(T^{-\lambda})^{-1}p_0$ is a horospherical net --- with hyperbolic Gauss map
$n$ --- if and only if $f^{\sharp}:\Z^2\to S^3\setminus(\C j\cup\{\infty\})$,
i.e.\ if $(T^{-\lambda})_{m,n}^{-1}p_0\not\in\C j\cup\{\infty\}$ at {\it one}
point $(m,n)$.  Clearly, {\it any} horospherical net can be (re-) constructed
from its hyperbolic Gauss map this way:

\Theorem{}{Let ${\fk g},{\fk h}:\Z^2\to\C$ be a Christoffel pair of
discrete isothermic nets, let ${\fk p}_0\in{\rm Im}\H\setminus\C j$,
and $\tau^{\lambda}:\Z^2\to Gl(2,\H)$, $\tau^{\lambda}_{0,0}=id$, be
a solution of $$\matrix{
	\tau^{\lambda}_{m+1,n} &=& \tau^{\lambda}_{m,n}(1+\lambda
		[v_0+v_{\infty}{\fk g}_{m,n}]
		(\partial_1{\fk h})_{m,n}
		[\nu_0-{\fk g}_{m+1,n}\nu_{\infty}]), \cr
	\tau^{\lambda}_{m,n+1} &=& \tau^{\lambda}_{m,n}(1+\lambda
		[v_0+v_{\infty}{\fk g}_{m,n}]
		(\partial_2{\fk h})_{m,n}
		[\nu_0-{\fk g}_{m,n+1}\nu_{\infty}]). \cr
	}\eqno{\rm(H)}$$
Then, $f^{\sharp}=(v_{\infty}\nu_0-v_0j\nu_{\infty})(\tau^{-\lambda})^{-1}
(v_0+v_{\infty}{\fk p}_0):\Z^2\to S^3\setminus S^2$ is a horospherical
net with hyperbolic Gauss map $n=v_0+v_{\infty}{\fk g}j:\Z^2\to S^2
=\C j\cup\{\infty\}$.\par
Every horospherical net $f^{\sharp}:\Z^2\to S^3\setminus(\C j\cup\{\infty\})$
can be constructed this way.}

On the other hand, horospherical nets are $T$-transforms of minimal nets
in Euclidean 3-space ${\rm Im}\H$:
with the discussion above, on the discrete Weierstrass representation for
minimal nets, and the $(T^{\lambda}C=D_{-\lambda}T^{\lambda})$-permutability
theorem, we conclude that any solution $T^{\lambda}$ of \WSystem\ provides
a horospherical net via $f=T^{\lambda}(v_0+v_{\infty}i){1\over\sqrt{2}}$
--- recall that \WSystem\ (or, (T)) does not depend on the stereographic
projection of $n=v_0+v_{\infty}{\fk g}j$ used to define the $C$-transform.
Thus,

\Theorem{}{Let ${\fk g}:\Z^2\to\C$ be an isothermic net, $\tau^{\lambda}$,
$\lambda\neq0$,
a solution of (H) as above; then, $f=(v_{\infty}\nu_0-v_0j\nu_{\infty})
\tau^{\lambda}(v_{\infty}i+v_0j){1\over\sqrt{2}}$ is a horospherical net
with hyperbolic Gauss map $n^{\lambda}=(v_{\infty}\nu_0-v_0j\nu_{\infty})
\tau^{\lambda}(v_0+v_{\infty}{\fk g})j:\Z^2\to S^2=\C j\cup\{\infty\}$.}

Note, that up to a change of model for the hyperbolic space, this is the
{\it exact} discrete analog of Bryant's Weierstrass type representation
for cmc-1 surfaces in hyperbolic space \UmYa.

\Definition{}{We refer to ${\fk g}$ as the secondary Gauss map of the
horospherical net $f$, and to the minimal net
$[i(i+{\fk g}j)(i-{\fk g}j)^{-1}]^{\ast}$ as the minimal cousin of $f$.}

As an example, Figure \CatCous\ indicates the family of horospherical nets with
the discrete catenoid (middle picture in figure \CatCous) as minimal cousins,
obtained from the discrete version of Bryant's Weierstrass type representation:
the ``discrete holomorphic nets'' ${\fk g}_{m,n}:=e^{2\pi{m+in\over N}}$ and
${\fk h}_{m,n}=e^{-2\pi{m+in\over N}}$ form a Christoffel pair since, by a
straightforward calculation, the (constant) cross ratios of elementary
quadrilaterals of ${\fk g}$ satisfy $q=
{\partial_1{\fk h}\partial_1{\fk g}\over\partial_2{\fk h}\partial_2{\fk g}}$.

Since the ``difference'' between two $D_{\lambda}$-transforms of
an isothermic net is measured by a $G$-transformation, the minimal
cousins of two horospherical nets with the same hyperbolic Gauss map
are $G$-transforms of each other.

The $T$-transforms of an isothermic net $n:\Z^2\to S^2$ take values in $S^2$.
Hence, the $(T^{\mu}D_{\lambda}=D_{\lambda-\mu}T^{\mu})$-permutability
theorem shows that (regular) $T$-transforms of horospherical nets are
horospherical. Moreover, all regular $T$-transforms of a horospherical
net have the same minimal cousin by the ``1-parameter group property'' of
the $T$-transformation, $T^{\lambda_2}T^{\lambda_1}=T^{\lambda_1+\lambda_2}$.

In \HHP, we showed that discrete minimal nets allow $\infty^3$ $D$-transforms
into minimal nets; this result carries over for horospherical nets via the
$T$-transformation.

By the compatibility of the $C$-$D_{\lambda}$-permutability with the
$T^{\lambda}$-transformation (see figure \PermThm), the roles of the
secondary and hyperbolic Gauss maps are interchanged for the horospherical
nets $f^{\sharp}$ and $f$ obtained from a ``holomorphic net'' ${\fk g}$ by
the two different representations given above, $n^{\lambda}=n^{\sharp}$;
therefore, $f$ and $f^{\sharp}$ can be considered as ``dual'' horospherical
nets (cf.\UmYaD).


\bigskip\noindent{\bf Acknowledgements}
I would like to thank F.~Burstall, F.~Pedit and U.~Pinkall for their interest
in my work, and for discussing their results with me.
In particular, the fabulous system \TSystem\ for the $T$-transformation is due
to discussions with F.~Burstall about the ``quaternionic function theory'' by
F.~Pedit and U.~Pinkall.


\References\bye

\Section{Appendix: Program text}
This is the program for producing the pictures in figure \CatCous\
using Mathematica:

{\petit\Verbatim

(* A graphics routine for quaternionic nets                   *)
(* ========================================================== *)

SURFACE[A_] := Block[ {gr={}} ,
   Do[ AppendTo[gr, Polygon[ {A[[i,j]][[2]],
                              A[[i+1,j]][[2]],
                              A[[i+1,j+1]][[2]],
                              A[[i,j+1]][[2]]} ] ],
   {i,1,Dimensions[A][[1]]-1} , {j,1,Dimensions[A][[2]]-1} ];
   Graphics3D[gr]
];

(* The program takes the parameter l,                         *)
(* step size k and ranges irg, jrg                            *)
(* ========================================================== *)

CATCOUS[l_,k_,irg_,jrg_] := Block[ {g,h,X,Y,TAU,TTransform} ,

(* The C-pair of discrete holomorphic nets                    *)
(* ========================================================== *)

g = Table[ N[ Exp[ 2 Pi (j + I i)/k ] ] ,
        {i,-irg,irg} , {j,-jrg,jrg} ];
h = Table[ N[ Exp[-2 Pi (j + I i)/k ] ] ,
        {i,-irg,irg} , {j,-jrg,jrg} ];

(* checked for isothermic dual nets: is ok ...
checkisoth = Table[ Chop[ N[
  (g[[i,j]]-g[[i+1,j]])/(g[[i+1,j]]-g[[i+1,j+1]]) *
  (g[[i+1,j+1]]-g[[i,j+1]])/(g[[i,j+1]]-g[[i,j]]) -
  ((g[[i+1,j]]-g[[i,j]]) (h[[i+1,j]]-h[[i,j]])) /
  ((g[[i,j+1]]-g[[i,j]]) (h[[i,j+1]]-h[[i,j]])) ] ] ,
  {i,1,2 irg} , {j,1,2 jrg} ];
*)

(* The connection form matrices                               *)
(* ========================================================== *)

X = Table[ N[
  {{1,0},{0,1}} + l*
  {{ g[[i,j]] (h[[i+1,j]]-h[[i,j]]) ,
    -g[[i,j]] (h[[i+1,j]]-h[[i,j]]) g[[i+1,j]] }, {
              (h[[i+1,j]]-h[[i,j]]) ,
             -(h[[i+1,j]]-h[[i,j]]) g[[i+1,j]] }} ] ,
  {i,1,2 irg} , {j,1,2 jrg + 1} ];

Y = Table[ N[
  {{1,0},{0,1}} + l*
  {{ g[[i,j]] (h[[i,j+1]]-h[[i,j]]) ,
    -g[[i,j]] (h[[i,j+1]]-h[[i,j]]) g[[i,j+1]] } , {
              (h[[i,j+1]]-h[[i,j]]) ,
             -(h[[i,j+1]]-h[[i,j]]) g[[i,j+1]] }} ] ,
  {i,1,2 irg + 1} , {j,1,2 jrg} ];

(* Integration of the system (H)                              *)
(* ========================================================== *)

TAU = Block[ {sol} ,
  sol = Table[ {{1,0},{0,1}} , {i,-irg,irg} , {j,-jrg,jrg} ];
  For[ i=irg+1 , i<=2 irg , i++ ,
    sol[[i+1,jrg+1]] = sol[[i,jrg+1]].X[[i,jrg+1]] ];
  For[ i=irg+1 , i>=2 , i-- ,
    sol[[i-1,jrg+1]] = sol[[i,jrg+1]].Inverse[ X[[i-1,jrg+1]] ] ];
  For[ i=1 , i<=2 irg + 1 , i++ ,
  For[ j=jrg+1 , j<=2 jrg , j++ ,
    sol[[i,j+1]] = sol[[i,j]].Y[[i,j]] ];
  For[ j=jrg+1 , j>=2 , j-- ,
    sol[[i,j-1]] = sol[[i,j]].Inverse[ Y[[i,j-1]] ] ];
  ];
  sol
  ];

(* The T-tranforms of g (hyperbolic Gauss maps)               *)
(* ========================================================== *)

TTransform = Block[ {xx} ,
  xx = Table[ { 0, { 0,0,0 }} , {i,-irg,irg} , {j,-jrg,jrg} ];
  For[ i=1 , i<=2 irg + 1 , i++ ,
  For[ j=1 , j<=2 jrg + 1 , j++ ,
    xx[[i,j]] = N[  { 0 , { 0 ,
      Re[ ( TAU[[i,j]][[1,1]] g[[i,j]] + TAU[[i,j]][[1,2]] )/
          ( TAU[[i,j]][[2,1]] g[[i,j]] + TAU[[i,j]][[2,2]] ) ] ,
      Im[ ( TAU[[i,j]][[1,1]] g[[i,j]] + TAU[[i,j]][[1,2]] )/
          ( TAU[[i,j]][[2,1]] g[[i,j]] + TAU[[i,j]][[2,2]] ) ] }}
  ] ] ];
  xx
  ];

(* The Catenoid cousins                                       *)
(* ========================================================== *)

CCousin = Block[ {xx} ,
  xx = Table[ { 0, { 0,0,0 }} , {i,-irg,irg} , {j,-jrg,jrg} ];
  For[ i=1 , i<=2 irg + 1 , i++ ,
  For[ j=1 , j<=2 jrg + 1 , j++ ,
    xx[[i,j]] = N[ 
    {-Im[ Det[TAU[[i,j]]] ] ,
    { Re[ Det[TAU[[i,j]]] ] ,
      Re[ TAU[[i,j]][[1,1]] Conjugate[TAU[[i,j]][[2,1]]] +
          TAU[[i,j]][[1,2]] Conjugate[TAU[[i,j]][[2,2]]] ] ,
      Im[ TAU[[i,j]][[1,1]] Conjugate[TAU[[i,j]][[2,1]]] +
          TAU[[i,j]][[1,2]] Conjugate[TAU[[i,j]][[2,2]]] ] }} /
    (TAU[[i,j]][[2,1]] Conjugate[TAU[[i,j]][[2,1]]] +
     TAU[[i,j]][[2,2]] Conjugate[TAU[[i,j]][[2,2]]])
  ] ] ];
  xx
  ];

(* ...and their image in the Poincare ball model              *)
(* ========================================================== *)

Spheremodell = Block[ {xx} ,
  xx = Table[ { 0 , { 1 , 0 , 0 }} , {i,-irg,irg} , {j,-jrg,jrg} ];
  For[ i=1 , i<=2 irg + 1 , i++ ,
  For[ j=1 , j<=2 jrg + 1 , j++ ,
    xx[[i,j]] = (CCousin[[i,j]]+{0,{1,0,0}})/(
                (CCousin[[i,j]][[2,1]]+1)^2
                +CCousin[[i,j]][[2,2]]^2
                +CCousin[[i,j]][[2,3]]^2)
  ]];
  xx
  ];

(* ...and the program spills out three quaternionic nets      *)
(* ========================================================== *)
Chop[ {TTransform,CCousin,Spheremodell} ]
];

(* A routine to produce and save the pictures                 *)
(* ========================================================== *)

Pic[{l_,k_,irg_,jrg_},file_:"cc.xx"] := Block[ {X} ,
  X = CATCOUS[l,k,irg,jrg];
  Display[ file , GraphicsArray[ {
  Show[ SURFACE[ X[[1]] ] ,
    (*PlotRange->All,*)Boxed->False,AspectRatio->Automatic ],
  Show[ SURFACE[ X[[2]] ] ,
    PlotRange->All,Boxed->False,AspectRatio->Automatic ],
  Show[ SURFACE[ X[[3]] ] ,
    PlotRange->All,Boxed->False,AspectRatio->Automatic ]
  } ] ]
  ];

(* The pictures                                               *)
(* ========================================================== *)

Pic[ {-.8000000 , 20 , 10 , 10 } , catcous1 ]
Pic[ {-.1170000 , 20 , 10 , 10 } , catcous2 ]
Pic[ {-.0500000 , 20 , 10 , 10 } , catcous3 ]
Pic[ {-.0250000 , 20 , 10 , 10 } , catcous4 ]
Pic[ {0.0000001 , 20 , 10 , 10 } , catenoid ]
Pic[ {0.0100000 , 20 , 10 , 10 } , catcous5 ]
Pic[ {0.0250000 , 20 , 10 , 10 } , catcous6 ]
Pic[ {0.0850000 , 20 , 10 , 10 } , catcous7 ]
Pic[ {0.2500000 , 20 , 10 , 10 } , catcous8 ]

\EndVerbatim}

\References
\bye